\documentclass{article}
\usepackage{mathrsfs, amsmath, amsthm,amssymb,bm}% packages for math
\usepackage{graphicx, xcolor, tikz}% packages for figures
\usetikzlibrary{arrows.meta, intersections, calc}
\usepackage{tkz-graph, tkz-berge}% If you use the xcolor package, load that package before tkz-graph to avoid package conflicts.
\usepackage{caption, subcaption}
\usepackage{array}
\usepackage{cases}
\makeatletter
\def\th@plain{%
  \upshape %\itshape % body font
}
\makeatother

\makeatletter
\renewenvironment{proof}[1][\proofname]{\par
  \pushQED{\qed}%
  \normalfont \topsep6\p@\@plus6\p@\relax
  \trivlist
  \item[\hskip\labelsep
        \bfseries
    #1\@addpunct{.}]\ignorespaces
}{%
  \popQED\endtrivlist\@endpefalse
}
\makeatother

\newtheorem{theorem}{Theorem}

\newtheorem*{conjecture*}{Conjecture}
\newtheorem{case}{Case}
\newtheorem{subcase}{Subcase}[case]

\theoremstyle{definition}

\newtheorem{remark}{Remark}

\usepackage[left=25mm,top=25mm,bottom=25mm,right=25mm]{geometry}
\usepackage[T1]{fontenc} % https://tex.stackexchange.com/questions/664/why-should-i-use-usepackaget1fontenc
\usepackage{newtxtext}
\usepackage[varvw]{newtxmath}

\usepackage{paralist, tablists}% put these packages ahead of enumitem
\usepackage[inline]{enumitem}
\usepackage[square, numbers, sort&compress]{natbib}

\usepackage[pdftex,%
  bookmarks=true,%
  bookmarksnumbered=true, %% 把章节的标号放入书签中，缺省为flase
  bookmarksopen=true, %% 打开书签树
  plainpages=false,%
  pdfpagelabels,%
  colorlinks=true, %% 彩色文字链接打开，去掉方框
  linkcolor=blue, %% 目录链接的颜色为蓝色
  citecolor=blue,%
  anchorcolor=green,
  urlcolor= blue,
  breaklinks=true,
  hyperindex=true]{hyperref}%% 不能和cite宏包同时使用

\newcounter{Hcase}
\newcounter{Hclaim}

\newcommand{\resetcounter}{\stepcounter{Hcase}\setcounter{case}{0}\stepcounter{Hclaim}\setcounter{claim}{0}}

\def\int(#1){\mathrm{int}(#1)}
\def\ext(#1){\mathrm{ext}(#1)}
\def\Int(#1){\mathrm{Int}(#1)}
\def\Ext(#1){\mathrm{Ext}(#1)}
\def\ad(#1){\mathrm{ad}(#1)}
\def\mad(#1){\mathrm{mad}(#1)}
\def\la(#1){\mathrm{la}(#1)}

\begin{document}
\title{Minor stars in plane graphs with minimum degree five}
\author{Yangfan Li { }\quad Mengjiao Rao { }\quad Tao Wang\footnote{{\tt Corresponding
author: wangtao@henu.edu.cn; iwangtao8@gmail.com} } \\
{\small Institute of Applied Mathematics}\\
{\small Henan University, Kaifeng, 475004, P. R. China}}
\date{}
\maketitle
\begin{abstract}
The weight of a subgraph $H$ in $G$ is the sum of the degrees in $G$ of vertices of $H$. The {\em height} of a subgraph $H$ in $G$ is the maximum degree of vertices of $H$ in $G$. A star in a given graph is minor if its center has degree at most five in the given graph. Lebesgue (1940) gave an approximate description of minor $5$-stars in the class of normal plane maps with minimum degree five. In this paper, we give two descriptions of minor $5$-stars in plane graphs with minimum degree five. By these descriptions, we can extend several results and give some new results on the weight and height for some special plane graphs with minimum degree five.  
\end{abstract}

\section{Introduction}
A {\em normal plane map} (NPM for short) is a connected plane pseudograph in which loops and multiple edges are allowed, but the degree of each vertex and face is at least three. A {\em $3$-polytope} is a $3$-connected plane graph. Clearly, each $3$-polytope is a normal plane map. The class of normal plane maps with minimum degree five is denoted by $\mathbf{M_{5}}$, and the class of $3$-polytopes with minimum degree at least five is denoted by $\mathbf{P_{5}}$. A $(k_{1}, k_{2}, k_{3}, k_{4}, k_{5})$-star is a star with $\deg(v_{i}) \leq k_{i}$, where $v_{i}$s are neighbors of the center in any order. A $k$-star is a star with $k$ rays. A star in a given graph is {\em minor} if its center has degree at most five in the given graph. 

The {\em weight} of a subgraph $H$ in $G$ is the sum of $\deg_{G}(v)$ by taking over all $v\in V(H)$. The {\em height} of a subgraph $H$ in $G$ is the maximum degree of vertices of $H$ in $G$. Let $\Omega_{\Delta}$ be the minimum integer such that there is a minor $5$-star with weight at most $\Omega_{\Delta}$ in every plane graph with minimum degree five and maximum degree $\Delta$. 

In 1904, Wernicke \cite{MR1511242} proved that every $M_{5}$ has a $5$-vertex adjacent to a $6^{-}$-vertex, that is a $(5, 6)$-edge. This was strengthened by Franklin \cite{MR1506473} in 1922 to the existence of a minor $(6, 6)$-star, that is a $(6, 5, 6)$-path. In 1996, Jendrol' and Madaras \cite{MR1446358} gave a precise description of minor $3$-stars in $\mathbf{M_{5}}$: there is a minor $(6, 6, 6)$- or $(5, 6, 7)$-star. In 1998, Borodin and Woodall \cite{MR1687827} showed that the minimum weight of minor $4$-star in $\mathbf{M_{5}}$ is at most $30$, which is best possible. Furthermore, Borodin and Ivanova \cite{MR3061007} gave a tight description of minor $4$-stars in $\mathbf{M_{5}}$. 

In 1940, Lebesgue \cite{MR0001903} gave an approximate description of minor $5$-stars in $\mathbf{M_{5}}$, which implies that $\Omega_{\Delta} \leq \Delta + 31$, and $\Omega_{\Delta} \leq \Delta + 27$ for $\Delta \geq 41$. In 1998, Borodin and Woodall \cite{MR1687827} strengthened the first of these results to $\Omega_{\Delta} \leq \Delta + 30$. This result is sharp for $\Delta = 7$ due to Borodin \cite{MR1235299} and Jendrol'--Madaras \cite{MR1446358}, $\Delta = 9$ due to Borodin--Ivanova \cite{MR3061007}, $\Delta = 10$ due to Jendrol'--Madaras \cite{MR1446358}, $\Delta = 12$ due to Borodin--Woodall \cite{MR1687827}. Recently, Borodin and Ivanova \cite{MR3665077} showed that $\Omega_{8} = 38$, $\Omega_{11} = 41$ and $\Omega_{13} = 42$. Hence, Borodin--Woodall's bound $\Delta + 30$ is sharp for every integer $\Delta$ in $\{7, 8, \dots, 12\}$, and we have known the exact value of $\Omega_{\Delta}$ for every integer $\Delta$ with $7 \leq \Delta \leq 13$. On the other hand, it is known that $\Omega_{20} = 48$. 

As for Lebesgue's bound $\Omega_{\Delta} \leq \Delta + 27$ for $\Delta \geq 41$, it was strengthened by Borodin, Ivanova and Jensen \cite{MR3227046} to $\Omega_{\Delta} \leq \Delta + 27$ for $\Delta \geq 28$, and further by Borodin and Ivanova \cite{MR3548785} to $\Omega_{\Delta} \leq \Delta + 27$ for $\Delta \geq 24$. 

In general, the description of minor stars is unordered for the neighbors of the center. In this paper, we give two descriptions of neighbors of $5$-vertices in a cyclic order. A $\langle \kappa_{1}, \kappa_{2}, \kappa_{3}, \kappa_{4}, \kappa_{5} \rangle$-star is a star with center having degree five and the other vertices having degrees $\leq \kappa_{1}, \leq \kappa_{2}, \leq \kappa_{3}, \leq \kappa_{4}, \leq \kappa_{5}$ in a cyclic order. 

The first purpose of this paper is to give the following description of minor $5$-stars in plane graphs with minimum degree five. Furthermore, this description of minor $5$-stars can imply the tight description of minor $4$-stars. For a problem of complete (tight) description of minor $5$-stars in plane graphs with minimum degree five, we refer the reader to \cite{MR3061007}. 
\begin{theorem}\label{MR1}
If $G$ is a plane graph with minimum degree five, then $G$ contains at least one of the following minor $5$-stars: 
\begin{tabenum}[ ]
\item $\langle 5, 7, 7, 5, 17 \rangle$, \item $\langle 5, 7, 8, 5, 11 \rangle$, \item $\langle 5, 7, 5, 8, 8 \rangle$, \item $\langle 8, 5, 5, 11, 6 \rangle$, \item $\langle 8, 5, 5, 8, 7 \rangle$, \item $\langle 8, 5, 5, 7, 8 \rangle$, \par
\item $\langle 8, 5, 5, 6, 9 \rangle$, \item $\langle 5, 6, 5, 8, 11 \rangle$, \item $\langle 5, 6, 6, 5, \infty \rangle$, \item $\langle 5, 6, 6, 6, 17 \rangle$, \item $\langle 6, 6, 6, 6, 11 \rangle$, \item $\langle 6, 6, 6, 7, 8 \rangle$, \par
\item $\langle 6, 6, 7, 6, 8 \rangle$, \item $\langle 5, 6, 6, 11, 7 \rangle$, \item $\langle 5, 6, 11, 6, 7 \rangle$, \item $\langle 5, 6, 6, 8, 8 \rangle$, \item $\langle 5, 6, 8, 6, 8 \rangle$, \item $\langle 5, 7, 6, 8, 7 \rangle$, \par
\item $\langle 5, 6, 7, 7, 7 \rangle$, \item $\langle 5, 6, 6, 7, 11 \rangle$, \item $\langle 5, 6, 7, 6, 11 \rangle$, \item $\langle 5, 7, 6, 7, 8 \rangle$, \item $\langle 5, 7, 7, 6, 8 \rangle$, \item $\langle 5, 7, 6, 6, 14 \rangle$, \par
\item $\langle 5, 8, 6, 6, 11 \rangle$, \item $\langle 5, 5, 7, 6, 14 \rangle$, \item $\langle 5, 6, 7, 5, 35 \rangle$, \item $\langle 5, 6, 8, 5, 15 \rangle$, \item $\langle 5, 6, 9, 5, 10 \rangle$, \item $\langle 5, 6, 11, 5, 9 \rangle$, \par
\item $\langle 5, 5, 10, 5, 12 \rangle$, \item $\langle 5, 7, 11, 5, 8 \rangle$, \item $\langle 5, 6, 5, 7, 14 \rangle$, \item $\langle 5, 5, 9, 5, 17 \rangle$. \item\item\qed 
\end{tabenum}
\end{theorem}

The following theorems are immediate consequences of \autoref{MR1}. Recall that the bounds in \autoref{T2}, \ref{T3}, \ref{T5}--\ref{T7} are sharp.  
\begin{theorem}[Borodin and Ivanova \cite{MR3665077}]\label{T2}
Let $\Delta$ be an integer with $\Delta \geq 13$. Every $3$-polytope with minimum degree five and maximum degree $\Delta$ has a minor $5$-star with weight at most $\Delta + 29$. 
\end{theorem}

\begin{theorem}[Borodin and Woodall \cite{MR1687827}]\label{T3}
Every plane graph with minimum degree five has a minor $4$-star with weight at most $30$. 
\end{theorem}

\begin{theorem}[Jendrol' and Madaras \cite{MR1446358}]\label{T4}
Every plane graph with minimum degree five has a minor $(10, 10, 10, 10)$-star.
\end{theorem}

\begin{theorem}[Jendrol' and Madaras \cite{MR1446358}]\label{T5}
Every plane graph with minimum degree five has a minor $(5, 6, 7)$-star or a minor $(6, 6, 6)$-star.
\end{theorem}

\begin{theorem}[Franklin \cite{MR1506473}]\label{T6}
Every plane graph with minimum degree five has a $(6, 5, 6)$-path. 
\end{theorem}

\begin{theorem}[Wernicke \cite{MR1511242}]\label{T7}
Every plane graph with minimum degree five has a $(5, 6)$-edge. 
\end{theorem}

\begin{theorem}
Every plane graph with minimum degree five having no $7$-vertices and no $8$-vertices contains at least one of the following minor 5-stars: a $\langle 5, 6, 6, 5, \infty \rangle$-star, a $\langle 5, 6, 6, 6, 17 \rangle$-star, a $\langle 6, 6, 6, 6, 11 \rangle$-star, a $\langle 5, 6, 9, 5, 10 \rangle$-star, a $\langle 5, 6, 11, 5, 9 \rangle$-star, a $\langle 5, 5, 10, 5, 12 \rangle$-star, a $\langle 5, 5, 9, 5, 17 \rangle$-star. 
\end{theorem}

\begin{theorem}
Every plane graph with minimum degree five having no vertices of degrees from $7$ to $9$ contains at least one of the following minor 5-stars: a $\langle 5, 6, 6, 5, \infty \rangle$-star, a $\langle 5, 6, 6, 6, 17 \rangle$-star, a $\langle 6, 6, 6, 6, 11 \rangle$-star, a $\langle 5, 5, 10, 5, 12 \rangle$-star. 
\end{theorem}

\begin{theorem}
Every plane graph with minimum degree five having no vertices of degrees from $7$ to $11$ contains at least one of the following minor 5-stars: a $\langle 5, 6, 6, 5, \infty \rangle$-star, a $\langle 5, 6, 6, 6, 17 \rangle$-star, a $\langle 6, 6, 6, 6, 6 \rangle$-star. 
\end{theorem}

\begin{theorem}[Borodin and Ivanova \cite{MR3528631}]\label{NO6}
Every $3$-polytope with minimum degree five having neither vertices of degrees from $6$ to $9$ nor $\langle 5, 5, 5, 5, \infty \rangle$-star has a minor $5$-star of weight at most $42$ and a minor $5$-star of height at most $12$. 
\end{theorem}

\begin{remark}
A tight description "a $(5, 6, 6, 5, \infty)$-star, a $(5, 6, 6, 6, 15)$-star, a $(6, 6, 6, 6, 6)$-star" for $3$-polytope with minimum degree $5$ and without vertices of degrees from $7$ to $11$ was obtained in \cite{MR3811956}. By \autoref{MR1}, the desired minor $5$-star to achieve the upper bounds in \autoref{NO6} is a $\langle 5, 5, 10, 5, 12\rangle$-star. Borodin and Ivanova \cite{MR3528631} have provided a minor $5$-star to show the sharpness of the bounds in \autoref{NO6}. 
\end{remark}

The second purpose of this paper is to give another description of minor $5$-stars plane graph with minimum degree five.  
\begin{theorem}\label{MR2}
If $G$ is a plane graph with minimum degree five, then $G$ contains at least one of the following minor $5$-stars: 
\begin{tabenum}[ ]
\item $\langle 5, 5, 5, 7, 17 \rangle$, \item $\langle 7, 5, 5, 7, 11 \rangle$, \item $\langle 7, 5, 5, 8, 9 \rangle$, \item $\langle 7, 5, 5, 9, 8 \rangle$, \item $\langle 7, 5, 5, 11, 7 \rangle$, \item $\langle 5, 5, 5, 8, 11 \rangle$, \par
\item $\langle 8, 5, 5, 9, 7 \rangle$, \item $\langle 8, 5, 5, 11, 6 \rangle$, \item $\langle 6, 6, 6, 6, 11 \rangle$, \item $\langle 6, 6, 6, 7, 9 \rangle$, \item $\langle 6, 6, 7, 6, 9 \rangle$, \item $\langle 6, 6, 7, 7, 7 \rangle$, \par
\item $\langle 6, 7, 6, 7, 7 \rangle$, \item $\langle 5, 6, 6, 8, 9 \rangle$, \item $\langle 5, 6, 8, 6, 9 \rangle$, \item $\langle 5, 6, 7, 7, 9 \rangle$, \item $\langle 5, 6, 6, 7, 11 \rangle$, \item $\langle 5, 6, 7, 6, 11 \rangle$, \par
\item $\langle 5, 7, 6, 7, 9 \rangle$, \item $\langle 5, 7, 7, 6, 9 \rangle$, \item $\langle 5, 6, 6, 6, 17 \rangle$, \item $\langle 5, 7, 6, 6, 11 \rangle$, \item $\langle 5, 8, 6, 6, 10 \rangle$, \item $\langle 5, 9, 6, 6, 9 \rangle$, \par
\item $\langle 5, 5, 9, 6, 9 \rangle$, \item $\langle 5, 6, 6, 5, \infty \rangle$, \item $\langle 5, 6, 7, 5, 23 \rangle$, \item $\langle 5, 6, 8, 5, 15 \rangle$, \item $\langle 5, 6, 9, 5, 14 \rangle$, \item $\langle 5, 9, 5, 6, 10 \rangle$, \par
\item $\langle 5, 8, 5, 6, 11 \rangle$, \item $\langle 5, 7, 5, 6, 17 \rangle$, \item $\langle 5, 7, 7, 5, 11 \rangle$, \item $\langle 5, 7, 8, 5, 9 \rangle$, \item $\langle 5, 8, 5, 7, 9 \rangle$, \item $\langle 5, 7, 5, 7, 11 \rangle$, \par
\item $\langle 5, 7, 5, 8, 10 \rangle$, \item $\langle 5, 7, 5, 9, 9 \rangle$, \item $\langle 5, 5, 10, 5, 14 \rangle$,  \item $\langle 5, 5, 11, 5, 13 \rangle$. \item \item \qed
\end{tabenum}
\end{theorem}

The following theorems are immediate consequences of \autoref{MR2}. 
\begin{theorem}
If $G$ is a plane graph with minimum degree $5$ and maximum degree $\Delta \geq 16$, then $G$ has a minor $5$-star with weight at most $\Delta + 28$. 
\end{theorem}

\begin{theorem}
Every plane graph with minimum degree five and no $\langle 5, 6, 6, 5, \infty \rangle$-stars and no $\langle 5, 6, 7, 5, 23 \rangle$-stars has a minor $5$-star of weight at most $45$ and height at most $17$. 
\end{theorem}

\begin{theorem}
Every plane graph with minimum degree five and no $\langle 5, 5, 5, 5, \infty \rangle$-stars and no $6$-vertices and  no $7$-vertices has a minor $5$-star of weight at most $44$ and height at most $15$. 
\end{theorem}

\begin{theorem}[Borodin and Ivanova \cite{MR3548785}]\label{WH}
Every plane graph with minimum degree five and no $\langle 5, 6, 6, 5, \infty \rangle$-star has a minor $5$-star of weight at most $51$ and height at most $23$. 
\end{theorem}

\begin{remark}
By \autoref{MR2}, the desired minor $5$-star to achieve the upper bounds in \autoref{WH} is a $\langle 5, 6, 7, 5, 23\rangle$-star.  
\end{remark}

Note that our results do not require the "3-connected" condition for the plane graphs with minimum degree five, so the class of graphs we considered is a little bit bigger than $\mathbf{P_{5}}$. For other results related to minor stars, we refer the reader to \cite{MR3061006, MR3431405, MR3601028, MR3634131, MR2346422, MR3762353,MR3811956,Borodin_2018,Borodin2018,Borodin2019+}, and more general results on the theory of light subgraphs in plane graphs has been presented recently in a nice survey \cite{MR3004475}. We use the classic discharging method and the ideas in \cite{MR3548785, MR3665077, MR3528631} to give a proof of \autoref{MR1} in \autoref{sec:2} and a proof of \autoref{MR2} in \autoref{sec:3}. 

Notions: In a plane graph with minimum degree five, let $v_{1}, v_{2}, \dots, v_{\kappa}$ be the neighbors of a $\kappa$-vertex $v$ in a cyclic order. If $v_{i}$ is a $5$-vertex and two of $v_{i-1}, v_{i+1}$ are $6^{+}$-vertices, then $v_{i}$ is called a {\em strong} $5$-neighbor of $v$; if $v_{i}$ is a $5$-vertex but it is not a strong $5$-neighbor of $v$, then $v_{i}$ is called a {\em non-strong} $5$-neighbor; if $v_{i-1}, v_{i}, v_{i+1}$ are all $5$-vertices, then $v_{i}$ is called a {\em weak} $5$-neighbor of $v$; if $v_{i-2}, v_{i-1}, v_{i}, v_{i+1}, v_{i+2}$ are all $5$-vertices, then $v_{i}$ is called a {\em twice-weak} $5$-neighbor of $v$. 

\section{Proof of \autoref{MR1}}\label{sec:2}
Let $G$ be a connected counterexample to \autoref{MR1} with maximum number of edges. 

\begin{enumerate}[label = ($\ast_{1}$)]
\item The graph $G$ is a triangulation. 
\end{enumerate}
\begin{proof}[Proof of ($\ast_{1}$)]
Suppose that $w_{1}, w_{2}, w_{3}, w_{4}$ are four consecutive vertices on the boundary of a $4^{+}$-face. Since $G$ is a simple graph, we have that $w_{1} \neq w_{3}$ and $w_{2} \neq w_{4}$. Note that $G$ is also a plane graph, thus we have that $w_{1}w_{3} \notin E(G)$ or $w_{2}w_{4} \notin E(G)$, otherwise the two lines representing $w_{1}w_{3}$ and $w_{2}w_{4}$ would cross each other outside the $4^{+}$-face. But an insertion of a diagonal $w_{1}w_{3}$ or $w_{2}w_{4}$ into the $4^{+}$-face would create a simple counterexample with more edges, which contradicts the assumption of $G$. 
\end{proof}

Euler's formula $|V| - |E| + |F| = 2$ for $G$ can be rewritten as the following: 
\begin{equation*}
\sum_{v\,\in\,V(G)}\big(\deg(v) - 6\big) + \sum_{f\,\in\,F(G)}\big(2\deg(f) - 6\big) = \sum_{v\,\in\,V(G)}\big(\deg(v) - 6\big) = -12.
\end{equation*}

Initially, we give every vertex $v$ an initial charge $\mu(v) = \deg(v) - 6$, and give every face $f$ an initial charge $\mu(f) = 2\deg(f) - 6$. Note that every face has an initial charge zero and every vertex has a nonnegative initial charge except the $5$-vertices. Next, we redistribute the charges among the vertices, preserving their sum,  such that the final charge $\mu'(v)$ of every vertex $v$ is nonnegative, which contradicts the fact that the sum of the initial charges is negative. 

Let 
\[
\alpha(\kappa) =
\begin{cases}
\frac{\kappa - 6}{\kappa}, & \text{if $\kappa \geq 8$ and $\kappa \notin \{11, 13, 14\}$;} \vspace{.1cm}\\
\frac{2}{5}, & \text{if $\kappa = 11$;} \vspace{.1cm}\\
\frac{1}{2}, & \text{if $\kappa = 13, 14$.}
\end{cases}
\]

\subsection{Discharging rules}
\begin{enumerate}[label = \bf R1\alph*]
\item Each $7$-vertex sends $\frac{1}{3}$ to each strong $5$-neighbor.
\item Each $7$-vertex sends $\frac{1}{6}$ to each non-strong $5$-neighbor.
\end{enumerate}

\begin{enumerate}[label = \bf R2]
\item\label{R2} Let $ww_{1}w_{2}$ be a $3$-face with $\deg(w) = \kappa$, where $\kappa \geq 8$ and $\kappa \notin \{10, 11\}$. 
\begin{enumerate}[label = \bf R2\alph*]
\item If both $w_{1}$ and $w_{2}$ are $5$-vertices, then $w$ sends $\frac{\alpha(\kappa)}{2}$ to each of $w_{1}$ and $w_{2}$ through this face. 
\item If $w_{1}$ is a $5$-vertex and $w_{2}$ is a $6^{+}$-vertex, then $w$ sends $\alpha(\kappa)$ to $w_{1}$ through this face. 
\end{enumerate}
\end{enumerate}

\begin{enumerate}[label = \bf R3]
\item\label{R3} Let $w$ be a $\kappa$-vertex with $\kappa = 10, 11$. Each such vertex $w$ sends $\frac{2}{5}$ to each adjacent vertex. Let $w_{0}, w_{1}, w_{2}$ be three consecutive neighbors of $w$ in a cyclic order. Suppose that $w_{0}$ is a $6^{+}$-vertex and $w_{1}$ is a $5$-vertex. 
\begin{enumerate}[label = \bf R3\alph*]
\item If $w_{2}$ is a $6^{+}$-vertex, then $w_{0}$ transfers a charge of $\frac{1}{5}$ to $w_{1}$. 
\item If $w_{2}$ is a $5$-vertex, then $w_{0}$ transfers a charge of $\frac{1}{10}$ to each of $w_{1}$ and $w_{2}$. 
\end{enumerate}
\end{enumerate}

\begin{enumerate}[label = \bf R4]
\item\label{R4} Each $11$-vertex additionally sends $\frac{1}{10}$ to each twice-weak $5$-neighbor. 
\end{enumerate}

\begin{enumerate}[label = \bf R5]
\item\label{R5} Each $13$-vertex or $14$-vertex additionally sends $\frac{1}{20}$ to each weak $5$-neighbor. 
\end{enumerate}

\begin{enumerate}[label = \bf R6]
\item\label{BACK1} Suppose that $w$ is a $5$-vertex with neighbors $w_{0}, w_{1}, w_{2}, w_{3}, w_{4}$ in a cyclic order, and $w_{0}, w_{1}, w_{2}, w_{3}, w_{4}$ have degrees $\kappa_{0}, 5, \kappa_{2}, \kappa_{3}, 5$, respectively. 
\begin{enumerate}[label = \bf R6\alph*]
\item\label{BACK11} If $\kappa_{2}, \kappa_{3} \geq 8$ and $\kappa_{0} = 13$, then $w$ returns $\frac{1}{4}$ to $w_{0}$. 
\item\label{BACK12} If $\kappa_{2}, \kappa_{3} \geq 9$, $\kappa_{0} = 11$ and $w$ is a twice-weak $5$-neighbor of $w_{0}$, then $w$ returns $\frac{1}{2}$ to $w_{0}$.
\item\label{BACK13} If $\kappa_{2}, \kappa_{3} = 8$, $\kappa_{0} = 11$ and $w$ is a twice-weak $5$-neighbor of $w_{0}$, then $w$ returns $\frac{1}{4}$ to $w_{0}$. 
\item\label{BACK14} If $\kappa_{2}, \kappa_{3} \geq 9$ and $\kappa_{0} = 7$, then $w$ returns $\frac{1}{6}$ to $w_{0}$. 
\end{enumerate}
\end{enumerate}

\begin{enumerate}[label = \bf R7]
\item\label{R7} Suppose that $w$ is a $5$-vertex with neighbors $w_{1}, w_{2}, w_{3}, w_{4}, w_{5}$ in a cyclic order, and $w_{1}, w_{2}, w_{3}, w_{4}$ have degrees $\kappa_{1}, 5, 5, \kappa_{4}$, respectively. 
\begin{enumerate}[label = \bf R7\alph*]
\item\label{A1} If $\kappa_{1} \geq 12$ and $\kappa_{4} \geq 16$, then $w$ donates $\frac{1}{8}$ to $w_{3}$. 
\item\label{A2} If $\kappa_{1} \geq 12$ and $\kappa_{4} = 13, 14, 15$, then $w$ donates $\frac{1}{20}$ to $w_{3}$. 
\item\label{A3} If $\kappa_{1} = 8$ and $\kappa_{4} \leq 11$, then $w$ donates $\frac{1}{8}$ to $w_{2}$.   
\end{enumerate}
\end{enumerate}

\begin{remark}
By \ref{R3}, each $10$-vertex sends $\frac{4}{5}$ to each strong $5$-neighbor, sends $\frac{2}{5}$ to each adjacent twice-weak $5$-neighbor, and sends at least $\frac{1}{2}$ to any other $5$-neighbor. 
\end{remark}

\begin{remark}
By \ref{R3} and \ref{R4}, each $11$-vertex sends $\frac{4}{5}$ to each strong $5$-neighbor, and sends at least $\frac{1}{2}$ to any other $5$-neighbor. 
\end{remark}
\begin{figure}[h]
\renewcommand\thesubfigure{\bf R1a}
\subcaptionbox{}{\begin{tikzpicture}[scale =2]
\coordinate (O) at (0, 0);
\coordinate (B) at (-120:1);
\coordinate (C) at (0, -1);
\coordinate (D) at (-60:1);
\draw (O) node[above] {$7$}--(B) node[left] {$6^{+}$}--(C) node[below] {$5$}--(D) node[right] {$6^{+}$}--cycle;
\draw (0, 0)--(0, -1);
\draw[-Stealth] (0, 0) --(0, -0.7) node[right] {\small$\frac{1}{3}$};
\fill (O) circle (1pt)
(B) circle (1pt)
(C) circle (1pt)
(D) circle (1pt);
\end{tikzpicture}}
\renewcommand\thesubfigure{\bf R1b}
\subcaptionbox{}{\begin{tikzpicture}[scale =2]
\coordinate (O) at (0, 0);
\coordinate (B) at (-120:1);
\coordinate (C) at (0, -1);
\coordinate (D) at (-60:1);
\draw (O) node[above] {$7$}--(B) node[left] {$5$}--(C) node[below] {$5$}--(D)--cycle;
\draw (0, 0)--(0, -1);
\draw[-Stealth] (0, 0) --(0, -0.7) node[right] {\small$\frac{1}{6}$};
\fill (O) circle (1pt)
(B) circle (1pt)
(C) circle (1pt)
(D) circle (1pt);
\end{tikzpicture}}
\renewcommand\thesubfigure{\bf R2a}
\subcaptionbox{\label{MR1:subfig:R2a}}{\begin{tikzpicture}[scale =2]
\coordinate (O) at (0, 0);
\coordinate (B) at (-120:1);
\coordinate (C) at (-60:1);
\draw (O) node[above] {$\kappa$}--(B) node[below] {$5$}--(C) node[below] {$5$}--cycle;
\fill (O) circle (1pt)
(B) circle (1pt)
(C) circle (1pt);
\draw [-{Stealth [sep =4pt]}](O) to [bend left =30](B);
\node at ($(B) + (9pt, 3pt)$) {\scriptsize$\frac{\alpha(\kappa)}{2}$};
\draw [-{Stealth [sep =4pt]}](O) to [bend right =30](C);
\node at ($(C) + (-9pt, 3pt)$) {\scriptsize$\frac{\alpha(\kappa)}{2}$};
\end{tikzpicture}}
\renewcommand\thesubfigure{\bf R2b}
\subcaptionbox{\label{MR1:subfig:R2b}}{\begin{tikzpicture}[scale =2]
\coordinate (O) at (0, 0);
\coordinate (B) at (-120:1);
\coordinate (C) at (-60:1);
\draw (O) node[above] {$\kappa$}--(B) node[below] {$5$}--(C) node[below] {$6^{+}$}--cycle;
\fill (O) circle (1pt)
(B) circle (1pt)
(C) circle (1pt);
\draw [-{Stealth [sep =4pt]}](O) to [bend left =30](B);
\node at ($(B) + (9pt, 3pt)$) {\scriptsize$\alpha(\kappa)$};
\end{tikzpicture}}
\renewcommand\thesubfigure{\bf R3a}
\subcaptionbox{}{\begin{tikzpicture}[scale =2]
\coordinate (O) at (0, 0);
\coordinate (B) at (-120:1);
\coordinate (C) at (0, -1);
\coordinate (D) at (-60:1);
\draw (O) node[above] {$10, 11$}--(B) node[left] {$6^{+}$}--(C) node[below] {$5$}--(D) node[right] {$6^{+}$}--cycle;
\draw (0, 0)--(0, -1);
\fill (O) circle (1pt)
(B) circle (1pt)
(C) circle (1pt)
(D) circle (1pt);
\draw [-{Stealth [sep =4pt]}](B) [bend left =30] to (C);
\node at ($(B) + (12pt, 2pt)$) {\small$\frac{1}{5}$};
\end{tikzpicture}}
\renewcommand\thesubfigure{\bf R3b}
\subcaptionbox{}{\begin{tikzpicture}[scale =2]
\coordinate (O) at (0, 0);
\coordinate (B) at (-120:1);
\coordinate (C) at (0, -1);
\coordinate (D) at (-60:1);
\draw (O) node[above] {$10, 11$}--(B) node[left] {$6^{+}$}--(C) node[below] {$5$}--(D) node[right] {$5$}--cycle;
\draw (0, 0)--(0, -1);
\fill (O) circle (1pt)
(B) circle (1pt)
(C) circle (1pt)
(D) circle (1pt);
\draw [-{Stealth [sep =4pt]}](B) [bend left =30] to (C);
\node at ($(B) + (12pt, 2pt)$) {\scriptsize$\frac{1}{10}$};
\draw [-{Stealth [sep =4pt]}](B) [bend left =30] to (D);
\node at ($(D) + (-5pt, 2pt)$) {\scriptsize$\frac{1}{10}$};
\end{tikzpicture}}
\renewcommand\thesubfigure{\bf R4}
\subcaptionbox{}{\begin{tikzpicture}[scale =2]
\coordinate (O) at (0, 0);
\coordinate (B) at (-130:1);
\coordinate (C) at (-110:1);
\coordinate (W0) at (-90:1);
\coordinate (W1) at (-70:1);
\coordinate (W2) at (-50:1);
\draw (O) node[above] {$11$}--(B) node[below] {$5$}--(C) node[below] {$5$}--(W0) node[below] {$5$}--(W1) node[below] {$5$}--(W2) node[below] {$5$}--cycle;
\draw (O)--(W0);
\draw [-Stealth](O)--(0, -0.7) node[right] {\small$\frac{1}{10}$};
\draw (C)--(O)--(W1);
\fill (O) circle (1pt)
(B) circle (1pt)
(C) circle (1pt)
(W0) circle (1pt)
(W1) circle (1pt)
(W2) circle (1pt);
\end{tikzpicture}}
\renewcommand\thesubfigure{\bf R5}
\subcaptionbox{}{\begin{tikzpicture}[scale =2]
\coordinate (O) at (0, 0);
\coordinate (B) at (-120:1);
\coordinate (C) at (0, -1);
\coordinate (D) at (-60:1);
\draw (O) node[above] {$13, 14$}--(B) node[left] {$5$}--(C) node[below] {$5$}--(D) node[right] {$5$}--cycle;
\draw (0, 0)--(0, -1);
\draw[-Stealth] (0, 0) --(0, -0.7) node[right] {\small$\frac{1}{20}$};
\fill (O) circle (1pt)
(B) circle (1pt)
(C) circle (1pt)
(D) circle (1pt);
\end{tikzpicture}}
\renewcommand\thesubfigure{\bf R6a}
\subcaptionbox{}{\begin{tikzpicture}[scale =2]
\draw (90:0.55) node[above] {$13$}--(162:0.55) node[left] {$5$}--(234:0.55) node[below] {$8^{+}$}--(306:0.55) node[below] {$8^{+}$}--(18:0.55) node[right] {$5$}--cycle;
\coordinate (O) at (0, 0);
\foreach \x in {90,162,234,306,18}
{\draw (0,0)--(\x:0.55);
\fill (\x:0.55) circle (1pt);}
\fill (O) circle (1pt);
\draw[-Stealth] (0, 0) --(0, 0.3) node[right] {\small$\frac{1}{4}$};
\end{tikzpicture}}
\renewcommand\thesubfigure{\bf R6b}
\subcaptionbox{}{\begin{tikzpicture}[scale =2]
\coordinate (O) at (0, 0);
\coordinate (X) at (-130:1);
\coordinate (W1) at (-110:1);
\coordinate (W) at (-90:1);
\coordinate (W4) at (-70:1);
\coordinate (Y) at (-50:1);
\draw (O) node[above] {$11$}--(X) node[below] {$5$}--(W1) node[below] {$5$}--(W) node[below] {$5$}--(W4) node[below] {$5$}--(Y) node[below] {$5$}--cycle;
\draw [-Stealth] (0, -0.7) --(0, -0.35) node[right] {\small$\frac{1}{2}$};
\coordinate (W2) at (-100:1.3);
\coordinate (W3) at (-80:1.3);
\fill (O) circle (1pt)
(X) circle (1pt)
(Y) circle (1pt)
(W0) circle (1pt)
(W1) circle (1pt)
(W2) circle (1pt)
(W3) circle (1pt)
(W4) circle (1pt);
\draw (W1)--(W2)  node[left] {$9^{+}$}--(W3) node[right] {$9^{+}$}--(W4)--(O)--cycle;
\draw (W2)--(W)--(W3);
\draw (W)--(O);
\end{tikzpicture}}
\renewcommand\thesubfigure{\bf R6c}
\subcaptionbox{}{\begin{tikzpicture}[scale =2]
\coordinate (O) at (0, 0);
\coordinate (X) at (-130:1);
\coordinate (W1) at (-110:1);
\coordinate (W) at (-90:1);
\coordinate (W4) at (-70:1);
\coordinate (Y) at (-50:1);
\draw (O) node[above] {$11$}--(X) node[below] {$5$}--(W1) node[below] {$5$}--(W) node[below] {$5$}--(W4) node[below] {$5$}--(Y) node[below] {$5$}--cycle;
\draw [-Stealth] (0, -0.7) --(0, -0.35) node[right] {\small$\frac{1}{4}$};
\coordinate (W2) at (-100:1.3);
\coordinate (W3) at (-80:1.3);
\fill (O) circle (1pt)
(X) circle (1pt)
(Y) circle (1pt)
(W0) circle (1pt)
(W1) circle (1pt)
(W2) circle (1pt)
(W3) circle (1pt)
(W4) circle (1pt);
\draw (W1)--(W2)  node[left] {$8$}--(W3) node[right] {$8$}--(W4)--(O)--cycle;
\draw (W2)--(W)--(W3);
\draw (W)--(O);
\end{tikzpicture}}\renewcommand\thesubfigure{\bf R6d}
\subcaptionbox{}{\begin{tikzpicture}[scale =2]
\draw (90:0.55) node[above] {$7$}--(162:0.55) node[left] {$5$}--(234:0.55) node[below] {$9^{+}$}--(306:0.55) node[below] {$9^{+}$}--(18:0.55) node[right] {$5$}--cycle;
\coordinate (O) at (0, 0);
\foreach \x in {90,162,234,306,18}
{\draw (0,0)--(\x:0.55);
\fill (\x:0.55) circle (1pt);}
\fill (O) circle (1pt);
\draw[-Stealth] (0, 0) --(0, 0.3) node[right] {\small$\frac{1}{6}$};
\end{tikzpicture}}
\renewcommand\thesubfigure{\bf R7a}
\subcaptionbox{}{\begin{tikzpicture}[scale =2]
\draw (90:0.55) node[above] {}--(162:0.55) node[left] {$12^{+}$}--(234:0.55) node[below] {$5$}--(306:0.55) node[below] {$5$}--(18:0.55) node[right] {$16^{+}$}--cycle;
\coordinate (O) at (0, 0);
\foreach \x in {90,162,234,306,18}
{\draw (0,0)--(\x:0.55);
\fill (\x:0.55) circle (1pt);}
\fill (O) circle (1pt);
\draw[-Stealth] (0, 0) --(306:0.3) node[right] {\small$\frac{1}{8}$};
\end{tikzpicture}}\renewcommand\thesubfigure{\bf R7b}
\subcaptionbox{\label{MR1:subfig:R7b}}{\begin{tikzpicture}[scale =2]
\draw (90:0.55) node[above] {}--(162:0.55) node[left] {$12^{+}$}--(234:0.55) node[below] {$5$}--(306:0.55) node[below] {$5$}--(18:0.55) node[right] {$\kappa_{4}$}--cycle;
\coordinate (O) at (0, 0);
\foreach \x in {90,162,234,306,18}
{\draw (0,0)--(\x:0.55);
\fill (\x:0.55) circle (1pt);}
\fill (O) circle (1pt);
\draw[-Stealth] (0, 0) --(306:0.3) node[right] {\small$\frac{1}{20}$};
\end{tikzpicture}}
\renewcommand\thesubfigure{\bf R7c}
\subcaptionbox{}{\begin{tikzpicture}[scale =2]
\draw (90:0.55) node[above] {}--(162:0.55) node[left] {$8$}--(234:0.55) node[below] {$5$}--(306:0.55) node[below] {$5$}--(18:0.55) node[right] {$11^{-}$}--cycle;
\coordinate (O) at (0, 0);
\foreach \x in {90,162,234,306,18}
{\draw (0,0)--(\x:0.55);
\fill (\x:0.55) circle (1pt);}
\fill (O) circle (1pt);
\draw[-Stealth] (0, 0) --(234:0.3) node[right] {\small$\frac{1}{8}$};
\end{tikzpicture}}
\caption{Discharging rules for \autoref{MR1}. Note that $\kappa \geq 8$ and $\kappa \notin \{10, 11\}$ in \subref{MR1:subfig:R2a} and \subref{MR1:subfig:R2b}, while $\kappa_{4} \in \{13, 14, 15\}$ in \subref{MR1:subfig:R7b}.} 
\end{figure}
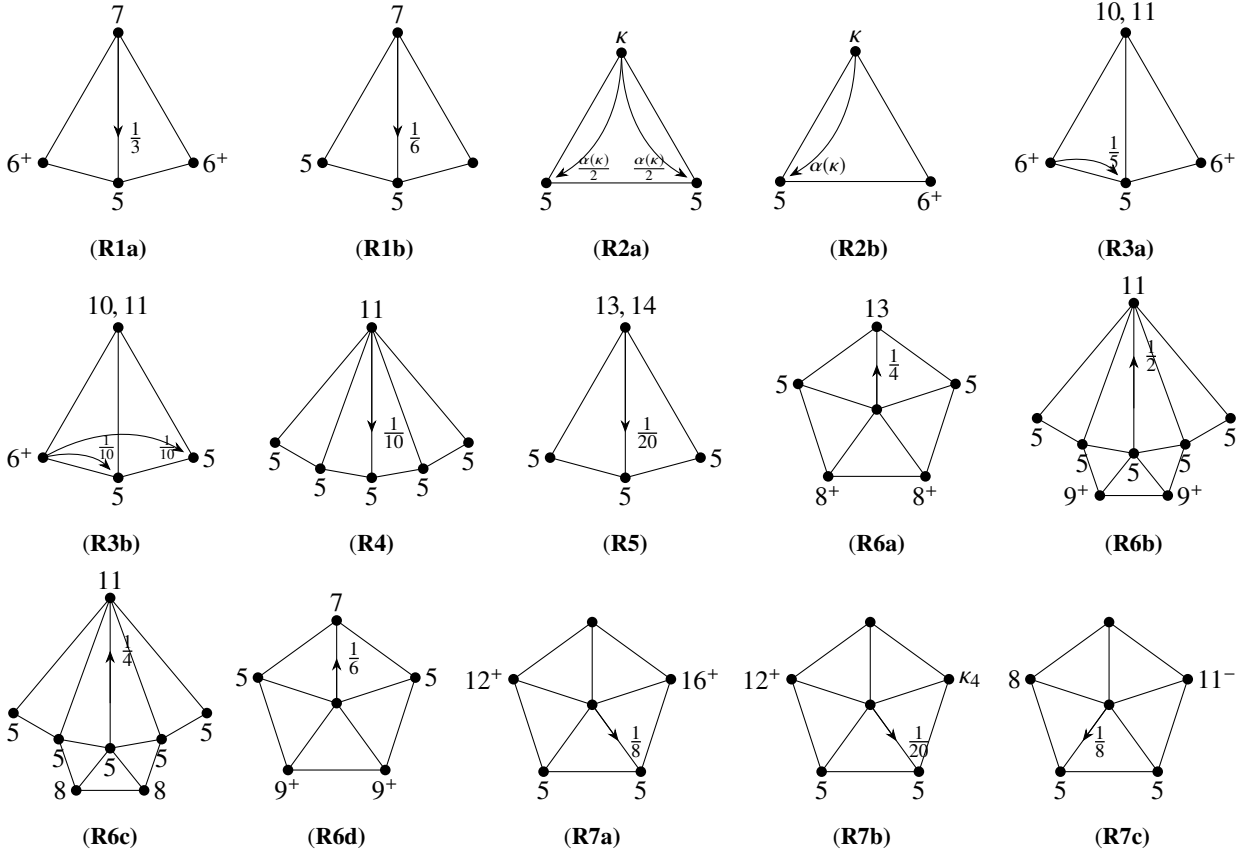
\subsection{The final charge of every vertex is nonnegative}
\begin{case}
If $v$ is a $\kappa$-vertex with $\kappa \geq 8$ and $\kappa \notin \{11, 13, 14\}$, then $\mu'(v) \geq \kappa - 6 - \kappa \cdot \alpha(\kappa) = 0$. 
\end{case}

\begin{case}
If $v$ is a $14$-vertex, then $\mu'(v) \geq 14 - 6 - 14 \cdot \alpha(14) - 14 \cdot \frac{1}{20} \geq 0$. 
\end{case}

\begin{case}
The vertex $v$ is a $13$-vertex. 
\end{case}
If $v$ has one $6^{+}$-neighbor, then it has at most ten weak $5$-neighbors, which implies that $\mu'(v) \geq 13 - 6 - 13 \cdot \alpha(13) - 10 \cdot \frac{1}{20} = 0$. So we may assume that $v$ is adjacent to thirteen $5$-vertices. By the oddness of the integer $13$ and the absence of $\langle 5, 7, 7, 5, 17 \rangle$-stars, the vertex $v$ is involved as the receiver in \ref{BACK11}, so $\mu'(v) \geq 13 - 6 - 13 \cdot (\frac{1}{2} + \frac{1}{20}) + \frac{1}{4} \geq 0$. 

\begin{case}
The vertex $v$ is an $11$-vertex. 
\end{case}
If $v$ has a $6^{+}$-neighbor, then it has at most six twice-weak neighbors, which implies that $\mu'(v) \geq 11 - 6 - 11 \cdot \frac{2}{5} - 6 \cdot \frac{1}{10} = 0$. It remains to assume that $v$ has eleven $5$-neighbors. If $v$ is involved as the receiver in \ref{BACK12}, then $\mu'(v) \geq 11 - 6 - 11 \cdot (\frac{2}{5} + \frac{1}{10}) + \frac{1}{2} = 0$. In the final case, the vertex $v$ is in a $\langle 5, 8, 8, 5, 11\rangle$-star due to the oddness of the integer $11$, so we have that $\mu'(v) \geq 11 - 6 - 11 \cdot (\frac{2}{5} + \frac{1}{10}) + 3 \cdot \frac{1}{4} \geq 0$ by \ref{BACK13}; otherwise, there is a $\langle 5, 7, 8, 5, 11 \rangle$-star. 

\begin{case}
The vertex $v$ is a $7$-vertex. 
\end{case}
If $v$ has at most three $5$-neighbors, then $\mu'(v) \geq 7 - 6 - 3 \cdot \frac{1}{3} = 0$. If $v$ has exactly four $5$-neighbors, then $v$ has exactly three $6^{+}$-neighbors and has at most two strong $5$-neighbors, which implies that $\mu'(v) \geq 7 - 6 - 2 \cdot \frac{1}{3} - 2 \cdot \frac{1}{6} = 0$. If $v$ has exactly five $5$-neighbors, then $v$ has at most one strong $5$-neighbor, which implies that $\mu'(v) \geq 7 - 6 - \frac{1}{3} - 4 \cdot \frac{1}{6} = 0$. If $v$ has exactly six $5$-neighbors, then $v$ has no strong $5$-neighbor, which implies that $\mu'(v) \geq 7 - 6 - 6 \cdot \frac{1}{6} = 0$. If $v$ has seven $5$-neighbors, then \ref{BACK14} is involved and $\mu'(v) \geq 7 - 6 - 7 \cdot \frac{1}{6} + \frac{1}{6} = 0$, for otherwise there is a $\langle 5, 7, 5, 8, 8 \rangle$-star. 

\begin{case}
The vertex $v$ is a $5$-vertex with neighbors $v_{1}, v_{2}, v_{3}, v_{4}, v_{5}$ in a cyclic order. Suppose that $v_{1}, v_{2}, v_{3}, v_{4}, v_{5}$ have degrees $\kappa_{1}, \kappa_{2}, \kappa_{3}, \kappa_{4}, \kappa_{5}$ respectively. Note that the initial charge of $v$ is $-1$. 
\end{case}

If $v$ is the sender in \ref{BACK11}, then $\mu'(v) \geq -1 + 2 \cdot \frac{3}{8} + \frac{1}{2} - \frac{1}{4} = 0$. If $v$ is the sender in \ref{BACK12}, then $\mu'(v) \geq -1 + 2 \cdot \frac{1}{2} + (\frac{2}{5} + \frac{1}{10}) - \frac{1}{2} = 0$. If $v$ is the sender in \ref{BACK13}, then $\mu'(v) \geq -1 + 2 \cdot \frac{3}{8} + (\frac{2}{5} + \frac{1}{10}) - \frac{1}{4} = 0$. If $v$ is the sender in \ref{BACK14}, then $\mu'(v) \geq -1 + 2 \cdot \frac{1}{2} + \frac{1}{6} - \frac{1}{6} = 0$. 

Suppose that $\kappa_{1}, \kappa_{4} \geq 12$ and $\kappa_{2} = \kappa_{3} = 5$. If $\kappa_{4} = 12$, then $v$ receives at least $\frac{1}{2}$ from $v_{4}$ and sends nothing to $v_{3}$. If $\kappa_{4} = 13, 14, 15$, then $v$ receives at least $\frac{11}{20}$ from $v_{4}$ and sends $\frac{1}{20}$ to $v_{3}$ by \ref{R2}, \ref{R5} and \ref{A2}, which implies that $v$ keeps at least $\frac{1}{2}$ from $v_{4}$. If $\kappa_{4} \geq 16$, then $v$ receives at least $\frac{5}{8}$ from $v_{4}$ and sends $\frac{1}{8}$ to $v_{3}$ by \ref{R2} and \ref{A1}, which implies that $v$ keeps at least $\frac{1}{2}$ from $v_{4}$. In conclusion, $v$ keeps at least $\frac{1}{2}$ from $v_{4}$ whenever $\kappa_{4} \geq 12$, and symmetrically $v$ keeps at least $\frac{1}{2}$ from $v_{1}$, which implies that $\mu'(v) \geq -1 + 2 \cdot \frac{1}{2} = 0$. 

Suppose that $\kappa_{1} = 8$, $\kappa_{2} = \kappa_{3} = 5$ and $\kappa_{4} \leq 11$. Then $v$ is involved as a sender in \ref{A3}. By the absence of $\langle 8, 5, 5, 11, 6 \rangle$-star, we have that $\kappa_{5} \geq 7$. If $\kappa_{5} = 7$, then $\kappa_{4} \geq 9$ and $\mu'(v) \geq -1 + (\frac{3}{8} - \frac{1}{8}) + \frac{1}{3} + \frac{1}{2} \geq 0$, for otherwise there is a $\langle 8, 5, 5, 8, 7 \rangle$-star. If $\kappa_{5} = 8$, then $\kappa_{4} \geq 8$ and $\mu'(v) \geq -1 + (\frac{3}{8} - \frac{1}{8}) + \frac{1}{2} + \frac{1}{4} \geq 0$, for otherwise there is a $\langle 8, 5, 5, 7, 8 \rangle$-star. If $\kappa_{5} = 9$, then $\kappa_{4} \geq 7$ and $\mu'(v) \geq -1 + (\frac{3}{8} - \frac{1}{8}) + \frac{2}{3} + \frac{1}{6} \geq 0$, for otherwise there is a $\langle 8, 5, 5, 6, 9 \rangle$-star. If $\kappa_{5} = 10, 11$, then $\kappa_{4} \geq 6$ and $\mu'(v) \geq -1 + (\frac{3}{8} - \frac{1}{8}) + \frac{4}{5} \geq 0$, for otherwise there is a $\langle 8, 5, 5, 5, 11 \rangle$-star. If $\kappa_{5} \geq 12 $, then $\mu'(v) \geq -1 + (\frac{3}{8} - \frac{1}{8}) + \frac{3}{4} \geq 0$.

So we may assume that the $5$-vertex $v$ is just a receiver in what follows. Suppose that $\kappa_{1}, \kappa_{2}, \kappa_{3}, \kappa_{4} \leq 6$, in other words, the vertex $v$ is the center of a $\langle 6, 6, 6, 6, \infty \rangle$-star. By the absence of $\langle 5, 6, 6, 5, \infty \rangle$-stars, we have that $\max\{\kappa_{1}, \kappa_{4}\} = 6$. If $\min\{\kappa_{1}, \kappa_{4}\} = 5$, then $\mu'(v) \geq -1 + \frac{12}{18} + \frac{12}{36}= -1 + \frac{3}{2} \cdot \frac{2}{3} = 0$, for otherwise there is a $\langle 5, 6, 6, 6, 17 \rangle$-star. If $\min\{\kappa_{1}, \kappa_{4}\} = 6$, then $\mu'(v) \geq -1 + 2 \cdot \frac{1}{2} = 0$, for otherwise there is a $\langle 6, 6, 6, 6, 11 \rangle$-star. So we may further assume that $v$ has at least two $7^{+}$-neighbors in what follows. 

\begin{subcase}
The $5$-vertex $v$ has no $5$-neighbor. 
\end{subcase} 

If $v$ has at least three $7^{+}$-neighbors, then $\mu'(v) \geq -1 + 3 \cdot \frac{1}{3} = 0$. If $v$ has at least two $8^{+}$-neighbors, then $\mu'(v) \geq -1 + 2 \cdot \frac{1}{2} = 0$. Hence, the vertex $v$ has exactly three $6$-neighbors and at most one $8^{+}$-neighbor, which implies that $\mu'(v) \geq -1 + \frac{1}{3} + \frac{2}{3} = 0$, for otherwise there is a $\langle 6, 6, 6, 7, 8 \rangle$-, or $\langle 6, 6, 7, 6, 8 \rangle$-star. 

\begin{subcase}
The $5$-vertex $v$ has precisely one $5$-neighbor $v_{1}$. 
\end{subcase}

By symmetry, we may assume that $\kappa_{3} \leq \kappa_{4}$. If $\kappa_{4} \geq 12$, then $\mu'(v) \geq -1 + 1 = 0$ and we are done. 

Suppose that $\kappa_{4} \in \{9, 10, 11\}$. If $\kappa_{3} \geq 7$, then $\mu'(v) \geq -1 + \frac{2}{3} + \frac{1}{3} = 0$. If $\min\{\kappa_{2}, \kappa_{5}\} \geq 7$, then $\mu'(v) \geq -1 + \frac{2}{3} + 2 \cdot \frac{1}{6} = 0$. If $\min\{\kappa_{2}, \kappa_{5}\} = \kappa_{3} = 6$, then $\mu'(v) \geq -1 + \frac{2}{3} + \frac{3}{8} \geq 0$, for otherwise there is a $\langle 5, 6, 6, 11, 7 \rangle$- or $\langle 5, 6, 11, 6, 7 \rangle$-star. 

Suppose that $v_{4}$ is an $8$-vertex. If $\kappa_{3} = 8$, then $\mu'(v) \geq -1 + 2 \cdot \frac{1}{2} = 0$. If $\kappa_{3} = 7$, then $\mu'(v) \geq -1 + \frac{1}{2} + \frac{1}{3} + \frac{1}{6} = 0$, for otherwise there is a $\langle 5, 6, 7, 8, 6 \rangle$-star. If $\min\{\kappa_{2}, \kappa_{5}\} = \kappa_{3} = 6$, then $\mu'(v) \geq -1 + 2 \cdot \frac{1}{2} = 0$, for otherwise there is a $\langle 5, 6, 6, 8, 8 \rangle$-, or $\langle 5, 6, 8, 6, 8 \rangle$-star. If $\kappa_{3} = 6$ and $\min\{\kappa_{2}, \kappa_{5}\} \geq 7$, then $\mu'(v) \geq -1 + \frac{1}{2} + \frac{1}{6} + \frac{3}{8} \geq 0$, for otherwise there is a $\langle 5, 7, 6, 8, 7 \rangle$-star. 

Suppose that $\kappa_{3} = \kappa_{4} = 7$. If $\min\{\kappa_{2}, \kappa_{5}\} \geq 7$, then $\mu'(v) \geq -1 + 2 \cdot \frac{1}{3} + 2 \cdot \frac{1}{6} = 0$. If $\min\{\kappa_{2}, \kappa_{5}\} = 6$, then $\mu'(v) \geq -1 + 2 \cdot \frac{1}{3} + \frac{3}{8} \geq 0$, for otherwise there is a $\langle 5, 6, 7, 7, 7 \rangle$-star. 

Suppose that $\kappa_{3} = 6$ and $\kappa_{4} = 7$. If $\min\{\kappa_{2}, \kappa_{5}\} = 6$, then $\mu'(v) \geq -1 + \frac{1}{3} + \frac{3}{4} \geq 0$, for otherwise there is a $\langle 5, 6, 6, 7, 11 \rangle$-, or $\langle 5, 6, 7, 6, 11 \rangle$-star. If $\min\{\kappa_{2}, \kappa_{5}\} = 7$, then $\mu'(v) \geq -1 + \frac{1}{3} + \frac{1}{6} + \frac{1}{2} = 0$, for otherwise there is a $\langle 5, 7, 6, 7, 8 \rangle$-, or $\langle 5, 7, 7, 6, 8 \rangle$-star. If $\min\{\kappa_{2}, \kappa_{5}\} \geq 8$, then $\mu'(v) \geq -1 + \frac{1}{3} + 2 \cdot \frac{3}{8} \geq 0$. 

Suppose that $\kappa_{3} = \kappa_{4} = 6$. If $\min\{\kappa_{2}, \kappa_{5}\} = 7$, then $\mu'(v) \geq -1 + \frac{1}{6} + \frac{3}{2} \cdot \frac{15 - 6}{15} \geq 0$, for otherwise there is a $\langle 5, 7, 6, 6, 14 \rangle$-star. If $\min\{\kappa_{2}, \kappa_{5}\} = 8$, then $\mu'(v) \geq -1 + \frac{3}{8} + \frac{3}{4} \geq 0$, for otherwise there is a $\langle 5, 8, 6, 6, 11 \rangle$-star. If $\min\{\kappa_{2}, \kappa_{5}\} \geq 9$, then $\mu'(v) \geq -1 + 2 \cdot \frac{1}{2} = 0$.

\begin{subcase}\label{C}
The $5$-vertex $v$ has precisely two $5$-neighbors $v_{2}$ and $v_{3}$.  
\end{subcase}

If $\min\{\kappa_{1}, \kappa_{4}\} \geq 9$, then $\mu'(v) \geq -1 + 2 \cdot \frac{1}{2} \geq 0$. If $\min\{\kappa_{1}, \kappa_{4}\} = 8$ and $\max\{\kappa_{1}, \kappa_{4}\} \geq 12$, then $\mu'(v) \geq -1 + \frac{3}{8} + \frac{3}{4} \geq 0$. Recall that the vertex $v$ is just a receiver, so we may assume that $\min\{\kappa_{1}, \kappa_{4}\} \leq 7$ and $\max\{\kappa_{1}, \kappa_{4}\} \neq 8$ in the following of Subcase~\ref{C}. 

Suppose that $\min\{\kappa_{1}, \kappa_{4}\} = 6$. If $\max\{\kappa_{1}, \kappa_{4}\} = 7$, then $\kappa_{5} \geq 12$ and $\mu'(v) \geq -1 + 1 + \frac{1}{6} \geq 0$, for otherwise there is a $\langle 6, 5, 5, 7, 11 \rangle$-star. If $\max\{\kappa_{1}, \kappa_{4}\} = 9, 10, 11$, then $\kappa_{5} \geq 8$ and $\mu'(v) \geq -1 + 2 \cdot \frac{1}{2} \geq 0$, for otherwise there is a $\langle 6, 5, 5, 11, 7 \rangle$-star. If $\max\{\kappa_{1}, \kappa_{4}\} \geq 12$, then $\mu'(v) \geq -1 + \frac{1}{3} + \frac{3}{4} \geq 0$. 

Suppose that $\min\{\kappa_{1}, \kappa_{4}\} = 7$. If $\max\{\kappa_{1}, \kappa_{4}\} = 7$, then $\kappa_{5} \geq 9$ and $\mu'(v) \geq -1 + 2 \cdot \frac{1}{6} + \frac{2}{3} \geq 0$, for otherwise there is a $\langle 7, 5, 5, 7, 8 \rangle$-star. If $9 \leq \max\{\kappa_{1}, \kappa_{4}\} \leq 14$, then $\kappa_{5} \geq 7$ and $\mu'(v) \geq -1 + \frac{1}{6} + \frac{1}{2} + \frac{1}{3} \geq 0$, for otherwise there is a $\langle 7, 5, 5, 14, 6 \rangle$-star. If $\max\{\kappa_{1}, \kappa_{4}\} \geq 15$, then $\mu'(v) \geq -1 + \frac{1}{6} + \frac{3}{2} \cdot \frac{3}{5} \geq 0$.

\begin{subcase}
The $5$-vertex $v$ has precisely two $5$-neighbors $v_{1}$ and $v_{3}$. 
\end{subcase}
As before, we may assume that $\kappa_{4} \leq \kappa_{5}$. If $\kappa_{4} \geq 9$, then $\mu'(v) \geq -1 + 2 \cdot \frac{1}{2} = 0$. 

Suppose that $v_{4}$ is a $6$-vertex. If $\kappa_{5} = 7$, then $\mu'(v) \geq -1 + \frac{1}{6} + \frac{36 - 6}{36} \geq 0$, for otherwise there is a $\langle 5, 6, 7, 5, 35 \rangle$-star. If $\kappa_{5} = 8$, then $\mu'(v) \geq -1 + \frac{3}{8} + \frac{16 - 6}{16} \geq 0$, for otherwise there is a $\langle 5, 6, 8, 5, 15 \rangle$-star. If $\kappa_{5} = 9$, then $\mu'(v) \geq -1 + 2 \cdot \frac{1}{2} = 0$, for otherwise there is a $\langle 5, 6, 9, 5, 10 \rangle$-star. If $\kappa_{5} \geq 15$, then $\mu'(v) \geq -1 + \frac{1}{6} + \frac{3}{2} \cdot \frac{15 - 6}{15} \geq 0$. If $12 \leq \kappa_{5} \leq 14$, then $\mu'(v) \geq -1 + \frac{3}{4} + \frac{1}{4} = 0$, for otherwise there is a $\langle 5, 7, 6, 6, 14 \rangle$-star. The remaining case is $\kappa_{5} = 10, 11$. By the absence of $\langle 5, 6, 11, 5, 9 \rangle$-stars, we have that $\kappa_{2} \geq 10$. If $v$ receives at least $\frac{1}{2}$ from $v_{2}$, then $\mu'(v) \geq -1 + 2 \cdot \frac{1}{2} = 0$. So we may assume that $\kappa_{2} = 10$ and $v$ is a twice-weak neighbor of $v_{2}$. Let $x, y, v_{2}, v, v_{5}$ be the neighbors of $v_{1}$ in a cyclic order. Since $v$ is a twice-weak neighbor of $v_{2}$, the vertex $y$ must be a $5$-vertex. Note that $v_{1}$ is not the center of a $\langle 5, 5, 10, 5, 12 \rangle$-star, thus $x$ cannot be a $5$-vertex. By \ref{R3}, we have that $\mu'(v) \geq -1 + \frac{2}{5} + (\frac{2}{5} + 2\cdot \frac{1}{10}) = 0$.

Suppose that $v_{4}$ is a $7$-vertex. If $\kappa_{5} = 7$, then $\mu'(v) \geq -1 + 2 \cdot \frac{1}{6} + \frac{18 - 6}{18} \geq 0$, for otherwise there is a $\langle 5, 7, 7, 5, 17 \rangle$-star. If $\kappa_{5} = 8$, then $\mu'(v) \geq -1 + \frac{1}{6} + \frac{3}{8} + \frac{1}{2} \geq 0$, for otherwise there is a $\langle 5, 7, 8, 5, 11 \rangle$-star. If $9 \leq \kappa_{5} \leq 11$, then $\mu'(v) \geq -1 + \frac{1}{6} + \frac{1}{2} + \frac{1}{3} = 0$, for otherwise there is a $\langle 5, 7, 11, 5, 8 \rangle$-star. If $12 \leq \kappa_{5} \leq 14$, then $\mu'(v) \geq -1 + 2 \cdot \frac{1}{6} + \frac{3}{4} \geq 0$, for otherwise there is a $\langle 5, 6, 5, 7, 14 \rangle$-star. If $\kappa_{5} \geq 15$, then $\mu'(v) \geq -1 + \frac{1}{6} + \frac{3}{2} \cdot \frac{15 - 6}{15} \geq 0$. 

Suppose that $v_{4}$ is an $8$-vertex. If $\kappa_{2} \geq 8$, then $\mu'(v) \geq -1 + 2 \cdot \frac{3}{8} + \frac{1}{4} = 0$. If $\kappa_{2} = 7$, then $\mu'(v) \geq -1 + \frac{1}{6} + \frac{3}{8} + \frac{1}{2} \geq 0$, for otherwise there is a $\langle 5, 7, 5, 8, 8 \rangle$-star. If $\kappa_{2} = 6$, then $\mu'(v) \geq -1 + \frac{3}{8} +  \frac{3}{4} \geq 0$, for otherwise there is a $\langle 5, 6, 5, 8, 11\rangle$-star.

\begin{subcase}
The $5$-vertex $v$ has precisely three $5$-neighbors $v_{1}, v_{2}$ and $v_{3}$. 
\end{subcase}

Note that $v$ is just a receiver and it has at least two $7^{+}$-neighbors, so we have that $\min\{\kappa_{4}, \kappa_{5}\} \neq 6, 8$. If $\min\{\kappa_{4}, \kappa_{5}\} = 7$, then $\mu'(v) \geq -1 + \frac{1}{6} + \frac{3}{2} \cdot \frac{15 - 6}{15} \geq 0$, for otherwise there is a $\langle 5, 6, 5, 7, 14 \rangle$-star. If $\min\{\kappa_{4}, \kappa_{5}\} \geq 9$, then $\mu'(v) \geq -1 + 2 \cdot \frac{1}{2} = 0$. 

\begin{subcase}
The $5$-vertex $v$ has precisely three $5$-neighbors $v_{2}, v_{3}$ and $v_{5}$.  
\end{subcase}

If $\min\{\kappa_{1}, \kappa_{4}\} \geq 11$, then $\mu'(v) \geq -1 + 2 \cdot \frac{1}{2} = 0$. If $\min\{\kappa_{1}, \kappa_{4}\} = 7$, then $\mu'(v) \geq -1 + \frac{1}{6} + \frac{36 - 6}{36} \geq 0$, for otherwise there is a $\langle 5, 6, 7, 5, 35 \rangle$-star. If $\min\{\kappa_{1}, \kappa_{4}\} = 9$, then $\mu'(v) \geq -1 + \frac{1}{3} + \frac{18 - 6}{18} \geq 0$, for otherwise there is a $\langle 5, 5, 9, 5, 17 \rangle$-star. 

\begin{figure}
\centering
\begin{tikzpicture}
\coordinate (V1) at (162:1);
\coordinate (V2) at (234:1);
\coordinate (V3) at (306:1);
\coordinate (V4) at (18:1);
\coordinate (V5) at (90:1);
\coordinate (W) at (-90:1.7);
\coordinate (X) at (220:2);
\coordinate (Y) at (-40:2);
\draw (V1) node[left]{\small$v_{1}$}--(V2) node[below]{\small$v_{2}$}--(V3) node[below]{\small$v_{3}$}--(V4) node[right]{\small$v_{4}$}--(V5) node[above]{\small$v_{5}$}--cycle;
\foreach \x in {162, 234, 306, 18, 90}
{\draw (0,0)--(\x:1);}
\draw (V1)--(X) node[left]{$x$}--(W) node[below]{$w$}--(Y) node[right]{$y$}--(V4);
\draw (X)--(V2)--(W)--(V3)--(Y);
\fill 
(V1) circle (2pt)
(V4) circle (2pt)
(W) circle (2pt)
(X) circle (2pt)
(Y) circle (2pt);
\fill [green](0, 0) circle (2pt)
(V2) circle (2pt)
(V3) circle (2pt)
(V5) circle (2pt);
\end{tikzpicture}
\caption{$v_{5}$ is a $5$-vertex}
\label{5585}
\end{figure}
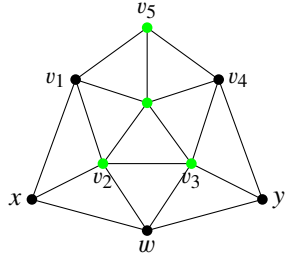

Suppose that $\min\{\kappa_{1}, \kappa_{4}\} = \kappa_{1} = 8$. Since there is no $\langle 5, 5, 8, 5, 15 \rangle$-stars, we must have that $\kappa_{4} \geq 16$, see \autoref{5585}. The $5$-vertex $v$ receives $\frac{1}{8}$ from $v_{2}$ whenever $w$ is a $11^{-}$-vertex by \ref{A1}, otherwise it receives $\frac{1}{8}$ from $v_{3}$ whenever $w$ is a $12^{+}$-vertex by \ref{A3}. Thus, $\mu'(v) \geq -1 + \frac{1}{4} + \frac{16 - 6}{16} + \frac{1}{8} = 0$. 

It suffices to consider $\min\{\kappa_{1}, \kappa_{4}\} = \kappa_{1} = 10$. By the absence of $\langle 5, 5, 10, 5, 12 \rangle$-stars, we have that $\kappa_{4} \geq 13$, see \autoref{5585}. If $v$ is not a twice-weak $5$-neighbor of $v_{1}$, then $\mu'(v) \geq -1 + (\frac{2}{5} + \frac{1}{10}) + \frac{1}{2} = 0$ by \ref{R2} and \ref{R3}. So we may assume that $v$ is a twice-weak $5$-neighbor of $v_{1}$, thus $\deg(w) \geq 13$ again by the absence of $\langle 5, 5, 10, 5, 12 \rangle$-stars. If $\kappa_{4} \in \{13, 14\}$, then $v$ receives $\frac{1}{2}$ from $v_{4}$ by \ref{R2}, and additionally $\frac{1}{20}$ from each of $v_{3}$ and $v_{4}$ by \ref{R5} and \ref{A2}, then $\mu'(v) \geq -1 + \frac{2}{5} + (\frac{1}{2} + \frac{1}{20}) + \frac{1}{20} = 0$. If $\kappa_{4} \geq 15$, then $v$ receives $\frac{2}{5}$ from $v_{1}$ and at least $\alpha(15) = \frac{3}{5}$ from $v_{4}$, then $\mu'(v) \geq -1 + \frac{2}{5} + \frac{3}{5} = 0$. This completes the proof of \autoref{MR1}. 
\resetcounter

\section{Proof of \autoref{MR2}}\label{sec:3}
Let $G$ be a connected counterexample to \autoref{MR2} with maximum number of edges. 

\begin{enumerate}[label = ($\ast_{2}$)]
\item The graph $G$ is a triangulation. 
\end{enumerate}
The proof of ($\ast_{2}$) is the same as that of ($\ast_{1}$), so we omit it. Euler's formula $|V| - |E| + |F| = 2$ for $G$ can be rewritten as the following: 
\begin{equation*}
\sum_{v\,\in\,V(G)}\big(\deg(v) - 6\big) + \sum_{f\,\in\,F(G)}\big(2\deg(f) - 6\big) = \sum_{v\,\in\,V(G)}\big(\deg(v) - 6\big) = -12.
\end{equation*}

Initially, we give every vertex $v$ an initial charge $\mu(v) = \deg(v) - 6$, and give every face $f$ an initial charge $\mu(f) = 2\deg(f) - 6$. Note that every face has an initial charge zero and every vertex has a nonnegative initial charge except the $5$-vertices. Next, we redistribute the charges among the vertices, preserving their sum,  such that the final charge $\mu'(v)$ of every vertex $v$ is nonnegative, which contradicts the fact that the sum of the initial charges is negative. 

Let
\[
\beta(\kappa) = 
\begin{array}{l@{\quad}r}
\frac{\kappa - 6}{\kappa},  &\text{if $\kappa \geq 8$.}
\end{array}
\]

\subsection{Discharging rules}
\begin{enumerate}[label = \bf R\arabic*]
\item Each $7$-vertex sends $\frac{1}{4}$ to each non-weak $5$-neighbor. 
\end{enumerate}

\begin{enumerate}[label = \bf R2]
\item Let $ww_{1}w_{2}$ be a $3$-face with $\deg(w) = \kappa$, where $\kappa \geq 8$ and $\kappa \neq 9$. 
\begin{enumerate}[label = \bf R2\alph*]
\item If both $w_{1}$ and $w_{2}$ are $5$-vertices, then $w$ sends $\frac{\beta(\kappa)}{2}$ to each of $w_{1}$ and $w_{2}$ through this face. 
\item If $w_{1}$ is a $5$-vertex and $w_{2}$ is a $6^{+}$-vertex, then $w$ sends $\beta(\kappa)$ to $w_{1}$ through this face. 
\end{enumerate}
\end{enumerate}

\begin{enumerate}[label = \bf R3]
\item\label{RR} Each $9$-vertex sends $\frac{1}{3}$ to each adjacent vertex. Let $w_{0}, w_{1}, w_{2}$ be three consecutive neighbors of a $9$-vertex in a cyclic order. Suppose that $w_{0}$ is a $6^{+}$-vertex and $w_{1}$ is a $5$-vertex.
\begin{enumerate}[label = \bf R3\alph*]
\item If $w_{2}$ is a $6^{+}$-vertex, then $w_{0}$ transfers a charge of $\frac{1}{6}$ to $w_{1}$. 
\item If $w_{2}$ is a $5$-vertex, then $w_{0}$ transfers a charge of $\frac{1}{12}$ to each of $w_{1}$ and $w_{2}$. 
\end{enumerate}
\end{enumerate}

\begin{enumerate}[label = \bf R4]
\item\label{BACK1} Suppose that $w$ is a $5$-vertex with neighbors $w_{1}, w_{2}, w_{3}, w_{4}, w_{5}$ in a cyclic order, and $w_{1}, w_{2}, w_{3}, w_{4}$ have degrees $\kappa_{1}, 5, 5, \kappa_{4}$, respectively. 
\begin{enumerate}[label = \bf R4\alph*]
\item\label{B1} If $\kappa_{1} \geq 12$ and $\kappa_{4} \geq 12$, then $w$ donates $\beta(\kappa_{1}) - \frac{1}{2}$ to $w_{2}$, and donates $\beta(\kappa_{4}) - \frac{1}{2}$ to $w_{3}$. 
\item\label{B2} If $\kappa_{1} = 7$ and $\kappa_{4} \leq 11$, then $w$ donates $\frac{1}{4}$ to $w_{2}$. 
\item\label{B3} If $\kappa_{1} = 8$ and $\kappa_{4} \leq 11$, then $w$ donates $\frac{1}{8}$ to $w_{2}$.   
\end{enumerate}
\end{enumerate}

\begin{remark}
By \ref{RR}, each $9$-vertex sends $\frac{2}{3}$ to each strong $5$-neighbor, sends $\frac{1}{3}$ to each twice-weak neighbor, and sends at least $\frac{5}{12}$ to any other $5$-neighbor. 
\end{remark}
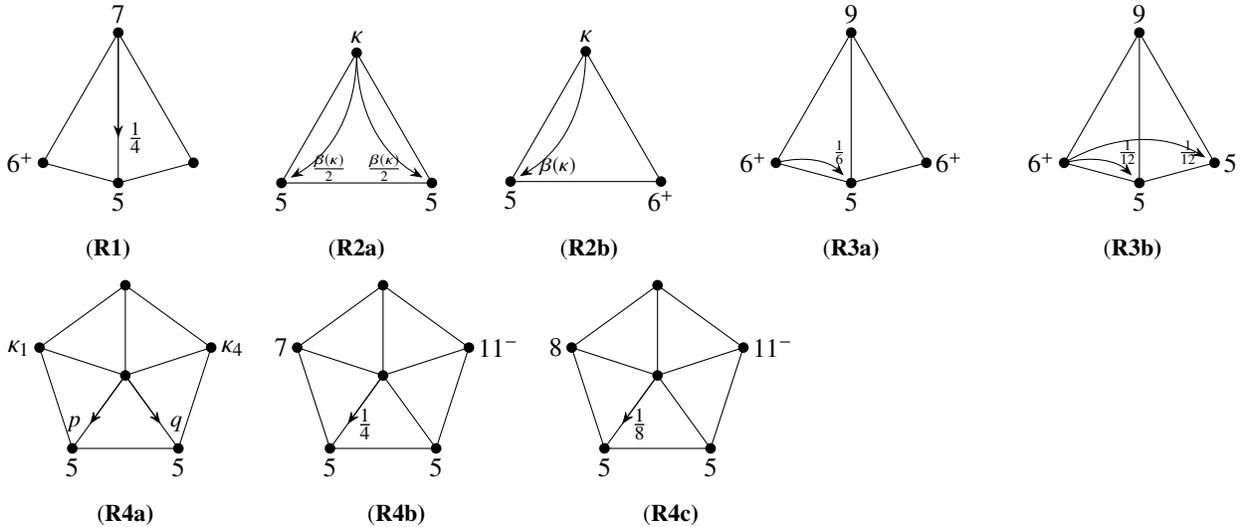
\begin{figure}[h]
\renewcommand\thesubfigure{\bf R1}
\subcaptionbox{}{\begin{tikzpicture}[scale =2]
\coordinate (O) at (0, 0);
\coordinate (B) at (-120:1);
\coordinate (C) at (0, -1);
\coordinate (D) at (-60:1);
\draw (O) node[above] {$7$}--(B) node[left] {$6^{+}$}--(C) node[below] {$5$}--(D) node[right] {$\mbox{ }$}--cycle;
\draw (0, 0)--(0, -1);
\draw[-Stealth] (0, 0) --(0, -0.7) node[right] {\small$\frac{1}{4}$};
\fill (O) circle (1pt)
(B) circle (1pt)
(C) circle (1pt)
(D) circle (1pt);
\end{tikzpicture}}
\renewcommand\thesubfigure{\bf R2a}
\subcaptionbox{\label{MR2:subfig:R2a}}{\begin{tikzpicture}[scale =2]
\coordinate (O) at (0, 0);
\coordinate (B) at (-120:1);
\coordinate (C) at (-60:1);
\draw (O) node[above] {$\kappa$}--(B) node[below] {$5$}--(C) node[below] {$5$}--cycle;
\fill (O) circle (1pt)
(B) circle (1pt)
(C) circle (1pt);
\draw [-{Stealth [sep =4pt]}](O) to [bend left =30](B);
\node at ($(B) + (9pt, 3pt)$) {\scriptsize$\frac{\beta(\kappa)}{2}$};
\draw [-{Stealth [sep =4pt]}](O) to [bend right =30](C);
\node at ($(C) + (-9pt, 3pt)$) {\scriptsize$\frac{\beta(\kappa)}{2}$};
\end{tikzpicture}}
\renewcommand\thesubfigure{\bf R2b}
\subcaptionbox{\label{MR2:subfig:R2b}}{\begin{tikzpicture}[scale =2]
\coordinate (O) at (0, 0);
\coordinate (B) at (-120:1);
\coordinate (C) at (-60:1);
\draw (O) node[above] {$\kappa$}--(B) node[below] {$5$}--(C) node[below] {$6^{+}$}--cycle;
\fill (O) circle (1pt)
(B) circle (1pt)
(C) circle (1pt);
\draw [-{Stealth [sep =4pt]}](O) to [bend left =30](B);
\node at ($(B) + (9pt, 3pt)$) {\scriptsize$\beta(\kappa)$};
\end{tikzpicture}}
\renewcommand\thesubfigure{\bf R3a}
\subcaptionbox{}{\begin{tikzpicture}[scale =2]
\coordinate (O) at (0, 0);
\coordinate (B) at (-120:1);
\coordinate (C) at (0, -1);
\coordinate (D) at (-60:1);
\draw (O) node[above] {$9$}--(B) node[left] {$6^{+}$}--(C) node[below] {$5$}--(D) node[right] {$6^{+}$}--cycle;
\draw (0, 0)--(0, -1);
\fill (O) circle (1pt)
(B) circle (1pt)
(C) circle (1pt)
(D) circle (1pt);
\draw [-{Stealth [sep =4pt]}](B) [bend left =30] to (C);
\node at ($(B) + (12pt, 2pt)$) {\scriptsize$\frac{1}{6}$};
\end{tikzpicture}}
\renewcommand\thesubfigure{\bf R3b}
\subcaptionbox{}{\begin{tikzpicture}[scale =2]
\coordinate (O) at (0, 0);
\coordinate (B) at (-120:1);
\coordinate (C) at (0, -1);
\coordinate (D) at (-60:1);
\draw (O) node[above] {$9$}--(B) node[left] {$6^{+}$}--(C) node[below] {$5$}--(D) node[right] {$5$}--cycle;
\draw (0, 0)--(0, -1);
\fill (O) circle (1pt)
(B) circle (1pt)
(C) circle (1pt)
(D) circle (1pt);
\draw [-{Stealth [sep =4pt]}](B) [bend left =30] to (C);
\node at ($(B) + (12pt, 2pt)$) {\scriptsize$\frac{1}{12}$};
\draw [-{Stealth [sep =4pt]}](B) [bend left =30] to (D);
\node at ($(D) + (-5pt, 2pt)$) {\scriptsize$\frac{1}{12}$};
\end{tikzpicture}}
\renewcommand\thesubfigure{\bf R4a}
\subcaptionbox{\label{MR2:subfig:R4a} }{\begin{tikzpicture}[scale =2]
\draw (90:0.6) node[above] {}--(162:0.6) node[left] {$\kappa_{1}$}--(234:0.6) node[below] {$5$}--(306:0.6) node[below] {$5$}--(18:0.6) node[right] {$\kappa_{4}$}--cycle;
\coordinate (O) at (0, 0);
\foreach \x in {90,162,234,306,18}
{\draw (0,0)--(\x:0.6);
\fill (\x:0.6) circle (1pt);}
\fill (O) circle (1pt);
\draw[-Stealth] (0, 0) --(234:0.4) node[left] {\small$p$};
\draw[-Stealth] (0, 0) --(306:0.4) node[right] {\small$q$};
\end{tikzpicture}}
\renewcommand\thesubfigure{\bf R4b}
\subcaptionbox{}{\begin{tikzpicture}[scale =2]
\draw (90:0.6) node[above] {}--(162:0.6) node[left] {$7$}--(234:0.6) node[below] {$5$}--(306:0.6) node[below] {$5$}--(18:0.6) node[right] {$11^{-}$}--cycle;
\coordinate (O) at (0, 0);
\foreach \x in {90,162,234,306,18}
{\draw (0,0)--(\x:0.6);
\fill (\x:0.6) circle (1pt);}
\fill (O) circle (1pt);
\draw[-Stealth] (0, 0) --(234:0.4) node[right] {\small$\frac{1}{4}$};
\end{tikzpicture}}
\renewcommand\thesubfigure{\bf R4c}
\subcaptionbox{}{\begin{tikzpicture}[scale =2]
\draw (90:0.6) node[above] {}--(162:0.6) node[left] {$8$}--(234:0.6) node[below] {$5$}--(306:0.6) node[below] {$5$}--(18:0.6) node[right] {$11^{-}$}--cycle;
\coordinate (O) at (0, 0);
\foreach \x in {90,162,234,306,18}
{\draw (0,0)--(\x:0.6);
\fill (\x:0.6) circle (1pt);}
\fill (O) circle (1pt);
\draw[-Stealth] (0, 0) --(234:0.4) node[right] {\small$\frac{1}{8}$};
\end{tikzpicture}}
\caption{Discharging rules for \autoref{MR2}. Note that $\kappa \geq 8$ but $\kappa \neq 9$ in \subref{MR2:subfig:R2a} and \subref{MR2:subfig:R2b}, while $\kappa_{1}, \kappa_{4} \geq 12$ and $p = \beta(\kappa_{1}) - \frac{1}{2}, q = \beta(\kappa_{4}) - \frac{1}{2}$ in \subref{MR2:subfig:R4a}.}
\end{figure}

\subsection{The final charge of every vertex is nonnegative}
\begin{case}
If $v$ is a $\kappa$-vertex with $\kappa \geq 8$ , then $\mu'(v) \geq \kappa - 6 - \kappa \cdot \beta(\kappa) = 0$. 
\end{case}

\begin{case}
The vertex $v$ is a $7$-vertex. 
\end{case}
If $v$ has at most four $5$-neighbors, then $\mu'(v) \geq 7 - 6 - 4 \cdot \frac{1}{4} = 0$. If $v$ has at least five $5$-neighbors, then it has at most two $6^{+}$-vertices and at most four non-weak $5$-neighbors, which also implies that $\mu'(v) \geq 7 - 6 - 4 \cdot \frac{1}{4} = 0$. 

\begin{case}
The vertex $v$ is a $5$-vertex with neighbors $v_{1}, v_{2}, v_{3}, v_{4}, v_{5}$ in a cyclic order. Suppose that $v_{1}, v_{2}, v_{3}, v_{4}$ and $v_{5}$ have degrees $\kappa_{1}, \kappa_{2}, \kappa_{3}, \kappa_{4}$ and $\kappa_{5}$ respectively. Note that the initial charge of $v$ is $-1$. 
\end{case}
If $v$ is the sender in \ref{B1} with $\kappa_{2} = \kappa_{3} = 5$, then $\mu'(v) \geq -1 + \beta(\kappa_{1}) + \beta(\kappa_{4}) - \big(\beta(\kappa_{1}) - \frac{1}{2}\big) - \big(\beta(\kappa_{4}) - \frac{1}{2}\big) = 0$. Suppose that $v$ is the sender in \ref{B2} with $\kappa_{2} = \kappa_{3} = 5$ and $\kappa_{1} = 7$. If $\kappa_{4} = 5$, then $\mu'(v) \geq -1 + \frac{1}{4} + \frac{3}{2} \cdot \frac{2}{3} - \frac{1}{4} = 0$, for otherwise there is a $\langle 5, 5, 5, 7, 17 \rangle$-star. If $\kappa_{4} = 6, 7$, then $\mu'(v) \geq -1 + \frac{1}{4} + 1 - \frac{1}{4} = -1 + 2 \cdot \frac{1}{4} + 1 - 2\cdot \frac{1}{4} = 0$, for otherwise there is a $\langle 7, 5, 5, 7, 11 \rangle$-star. If $\kappa_{4} = 8$, then $\mu'(v) \geq -1 + \frac{1}{4} + \frac{3}{8} + \frac{4}{5} - \frac{1}{4} - \frac{1}{8} \geq 0$, for otherwise there is a $\langle 7, 5, 5, 8, 9 \rangle$-star. If $\kappa_{4} = 9$, then $\mu'(v) \geq -1 + \frac{1}{4} + \frac{5}{12} + \frac{2}{3} - \frac{1}{4} \geq 0$, for otherwise there is a $\langle 7, 5, 5, 9, 8 \rangle$-star. If $\kappa_{4} = 10, 11$, then $\mu'(v) \geq -1 + \frac{1}{4} + \frac{3}{5} + \frac{1}{2} - \frac{1}{4} \geq 0$, for otherwise there is a $\langle 7, 5, 5, 11, 7 \rangle$-star.

Suppose that $v$ is the sender in \ref{B3} with $\kappa_{2} = \kappa_{3} = 5$ and $\kappa_{1} = 8$. If $\kappa_{4} = 5$, then $\mu'(v) \geq -1 + \frac{3}{8} + \frac{3}{2} \cdot \frac{1}{2} - \frac{1}{8} = 0$, for otherwise there is a $\langle 5, 5, 5, 8, 11 \rangle$-star. If $\kappa_{4} = 6$, then $\mu'(v) \geq -1 + \frac{3}{8} + \frac{4}{5} - \frac{1}{8} \geq 0$, for otherwise there is a $\langle 6, 5, 5, 8, 9 \rangle$-star. If $\kappa_{4} = 7$, then $\mu'(v) \geq -1 + \frac{3}{8} + \frac{4}{5} + \frac{1}{4} - \frac{1}{8} - \frac{1}{4} \geq 0$, for otherwise there is a $\langle 7, 5, 5, 8, 9 \rangle$-star. If $\kappa_{4} = 8$, then $\mu'(v) \geq -1 + 2 \cdot \frac{3}{8} + \frac{1}{2} - 2 \cdot \frac{1}{8} = 0$, for otherwise there is an $\langle 8, 5, 5, 8, 7 \rangle$-star. If $\kappa_{4} = 9$, then $\mu'(v) \geq -1 + \frac{3}{8} + \frac{5}{12} + \frac{1}{2} - \frac{1}{8} \geq 0$, for otherwise there is an $\langle 8, 5, 5, 9, 7 \rangle$-star. If $\kappa_{4} = 10, 11$, then $\mu'(v) \geq -1 + \frac{3}{8} + \frac{1}{4} + \frac{3}{5} - \frac{1}{8} \geq 0$, for otherwise there is an $\langle 8, 5, 5, 11, 6 \rangle$-star. 

So we may assume that the $5$-vertex $v$ is just a receiver in what follows.

\begin{subcase}
The $5$-vertex $v$ has no $5$-neighbor. 
\end{subcase} 

If $v$ has four $6$-neighbors, then $\mu'(v) \geq -1 + 1 = 0$, for otherwise there is a $\langle 6, 6, 6, 6, 11 \rangle$-star. If $v$ has at least two $8^{+}$-neighbors, then $\mu'(v) \geq -1 + 2 \cdot \frac{1}{2} = 0$. If $v$ has exactly three $6$-neighbors and one $7$-neighbor, then $\mu'(v) \geq -1 + \frac{1}{4} + \frac{4}{5} \geq 0$, for otherwise there is a $\langle 6, 6, 6, 7, 9 \rangle$-star or a $\langle 6, 6, 7, 6, 9 \rangle$-star. If $v$ has exactly two $6$-neighbors, then $\mu'(v) \geq -1 + 2 \cdot \frac{1}{4} + \frac{1}{2} = 0$, for otherwise there is a $\langle 6, 6, 7, 7, 7 \rangle$-star or a $\langle 6, 7, 6, 7, 7 \rangle$-star. If $v$ has at most one $6$-neighbor, then $\mu'(v) \geq -1 + 4 \cdot \frac{1}{4} = 0$. 

\begin{subcase}
The $5$-vertex $v$ has precisely one $5$-neighbor $v_{1}$. 
\end{subcase}

By symmetry, we may assume that $\kappa_{3} \leq \kappa_{4}$. If $\kappa_{3} \geq 8$, then $v$ receives at least $\frac{1}{2}$ from each of $v_{3}$ and $v_{4}$, which implies that $\mu'(v) \geq -1 + 2 \cdot \frac{1}{2} = 0$. If $\kappa_{4} \in \{10, 11\}$, then $\mu'(v) \geq -1 + \frac{4}{5} + \frac{1}{4} \geq 0$, for otherwise there is a $\langle 6, 6, 6, 6, 11 \rangle$-star. If $\kappa_{4} \geq 12$, then $\mu'(v) \geq -1 + 1 = 0$ and we are done. So we may assume that $\kappa_{3} \leq 7$ and $\kappa_{4} \leq 9$.

Suppose that $v_{4}$ is a $9$-vertex. If $\kappa_{3} = 7$, then $\mu'(v) \geq -1 + \frac{2}{3} + 2 \cdot \frac{1}{4} \geq 0$, for otherwise there is a $\langle 6, 5, 6, 7, 9 \rangle$-star. If $\min\{\kappa_{2}, \kappa_{5}\} = \kappa_{3} = 6$, then $\mu'(v) \geq -1 + \frac{2}{3} + \frac{3}{8} \geq 0$, for otherwise there is a $\langle 6, 6, 6, 7, 9 \rangle$-, or $\langle 6, 6, 7, 6, 9 \rangle$-star. If $\min\{\kappa_{2}, \kappa_{5}\} \geq 7$, then $\mu'(v) \geq -1 + \frac{2}{3} + 2 \cdot \frac{1}{4} \geq 0$. 

Suppose that $v_{4}$ is an $8$-vertex. If $\kappa_{3} = 7$, then $\mu'(v) \geq -1 + \frac{1}{2} + 2 \cdot \frac{1}{4} = 0$, for otherwise there is a $\langle 6, 6, 6, 7, 9 \rangle$-star. If $\kappa_{3} = 6$ and $\min\{\kappa_{2}, \kappa_{5}\} \geq 7$, then $\mu'(v) \geq -1 + \frac{1}{2} + 2 \cdot \frac{1}{4} = 0$. If $\min\{\kappa_{2}, \kappa_{5}\} = \kappa_{3} = 6$, then $\mu'(v) \geq -1 + \frac{1}{2} + \frac{3}{5} \geq 0$, for otherwise there is a $\langle 5, 6, 6, 8, 9 \rangle$- or $\langle 5, 6, 8, 6, 9 \rangle$-star.

Suppose that $\kappa_{3} = \kappa_{4} = 7$. If $\min\{\kappa_{2}, \kappa_{5}\} \geq 7$, then $\mu'(v) \geq -1 + 4 \cdot \frac{1}{4} = 0$. If $\min\{\kappa_{2}, \kappa_{5}\} = 6$, then $\mu'(v) \geq -1 + 2 \cdot \frac{1}{4} + \frac{3}{5} \geq 0$, for otherwise there is a $\langle 5, 6, 7, 7, 9 \rangle$-star. 

Suppose that $\kappa_{4} = 7$ and $\kappa_{3} = 6$. If $\min\{\kappa_{2}, \kappa_{5}\} = 6$, then $\mu'(v) \geq -1 + \frac{1}{4} + \frac{3}{4} = 0$, for otherwise there is a $\langle 5, 6, 6, 7, 11 \rangle$- or $\langle 5, 6, 7, 6, 11 \rangle$-star. If $\min\{\kappa_{2}, \kappa_{5}\} = 7$, then $\mu'(v) \geq -1 + 2 \cdot \frac{1}{4} + \frac{3}{5} \geq 0$, for otherwise there is a $\langle 5, 7, 6, 7, 9 \rangle$- or $\langle 5, 7, 7, 6, 9 \rangle$-star. If $\min\{\kappa_{2}, \kappa_{5}\} \geq 8$, then $\mu'(v) \geq -1 + 2 \cdot \frac{3}{8} + \frac{1}{4} = 0$. 

Suppose that $\kappa_{3} = \kappa_{4} = 6$. If $\min\{\kappa_{2}, \kappa_{5}\} = 6$, then $\mu'(v) \geq -1 + \frac{3}{2} \cdot \frac{18 - 6}{18} \geq 0$, for otherwise there is a $\langle 5, 6, 6, 6, 17 \rangle$-star. If $\min\{\kappa_{2}, \kappa_{5}\} = 7$, then $\mu'(v) \geq -1 + \frac{1}{4} + \frac{3}{4} = 0$, for otherwise there is a $\langle 5, 7, 6, 6, 11 \rangle$-star. If $\min\{\kappa_{2}, \kappa_{5}\} = 8$, then $\mu'(v) \geq -1 + \frac{3}{8} + \frac{15}{22} \geq 0$, for otherwise there is a $\langle 5, 8, 6, 6, 10 \rangle$-star. If $\min\{\kappa_{2}, \kappa_{5}\} \geq 9$, then $\mu'(v) \geq -1 + \frac{5}{12} + \frac{3}{5} \geq 0$, for otherwise there is a $\langle 5, 9, 6, 6, 9 \rangle$-star. 

\begin{subcase}\label{CC}
The $5$-vertex $v$ has precisely two $5$-neighbors $v_{2}$ and $v_{3}$.  
\end{subcase}

If $\min\{\kappa_{1}, \kappa_{4}\} \in \{7, 8\}$ and $\max\{\kappa_{1}, \kappa_{4}\} \geq 12$, then $\mu'(v) \geq -1 + \frac{1}{4} + \frac{3}{4} = 0$. Recall that the vertex $v$ is just a receiver, so we may assume that $\kappa_{1} \neq 7, 8$ and $\kappa_{4} \neq 7, 8$ in the following of Subcase~\ref{CC}. 

Suppose that $\min\{\kappa_{1}, \kappa_{4}\} = 6$. If $\kappa_{5} \geq 12$, then $\mu'(v) \geq -1 + 1 = 0$. If $\kappa_{5} = 10, 11$, then $\mu'(v) \geq -1 + \frac{4}{5} + \frac{1}{4} \geq 0$, for otherwise there is a $\langle 6, 5, 5, 6, 11 \rangle$-star. If $\kappa_{5} = 9$, then $\mu'(v) \geq -1 + \frac{2}{3} + \frac{3}{8} \geq 0$, for otherwise there is a $\langle 6, 5, 5, 7, 9\rangle$-star. If $\kappa_{5} = 8$, then $\mu'(v) \geq -1 + \frac{1}{2} + \frac{3}{5} \geq 0$, for otherwise there is a $\langle 5, 5, 6, 8, 9 \rangle$-star. If $\kappa_{5} = 7$, then $\mu'(v) \geq -1 + \frac{1}{4} + \frac{3}{4} = 0$, for otherwise there is a $\langle 5, 5, 6, 7, 11 \rangle$-star. If $\kappa_{5} = 6$, then $\mu'(v) \geq -1 + \frac{3}{2} \cdot \frac{2}{3} = 0$, for otherwise there is a $\langle 5, 5, 6, 6, 17 \rangle$-star. 

Suppose that $\min\{\kappa_{1}, \kappa_{4}\} \geq 9$. If $\max\{\kappa_{1}, \kappa_{4}\} \geq 10$, then $\mu'(v) \geq -1 + \frac{5}{12} + \frac{3}{5} \geq 0$. If $\max\{\kappa_{1}, \kappa_{4}\} = 9$, then $\kappa_{1} = \kappa_{4} = 9$, which implies that $\mu'(v) \geq -1 + 2 \cdot \frac{5}{12} + \frac{1}{4} \geq 0$, for otherwise there is a $\langle 5, 5, 9, 6, 9 \rangle$-star.

\begin{subcase}
The $5$-vertex $v$ has precisely two $5$-neighbors $v_{1}$ and $v_{3}$. 
\end{subcase}
As before, we may assume that $\kappa_{4} \leq \kappa_{5}$. If $\kappa_{4} \geq 10$, then $\mu'(v) \geq -1 + 2 \cdot \frac{3}{5} \geq 0$. 

Suppose that $v_{4}$ is a $6$-vertex. By the absence of $\langle 5, 6, 6, 5, \infty \rangle$-stars, we have that $\kappa_{5} \geq 7$. If $\kappa_{5} = 7$, then $\mu'(v) \geq -1 + \frac{1}{4} + \frac{24 - 6}{24} \geq 0$, for otherwise there is a $\langle 5, 6, 7, 5, 23 \rangle$-star. If $\kappa_{5} = 8$, then $\mu'(v) \geq -1 + \frac{3}{8} + \frac{5}{8} = 0$, for otherwise, there is a $\langle 5, 6, 8, 5, 15 \rangle$-star. If $\kappa_{5} = 9$, then $\mu'(v) \geq -1 + \frac{5}{12} + \frac{3}{5} \geq 0$, for otherwise there is a $\langle 5, 6, 9, 5, 14 \rangle$-star. If $\kappa_{5} = 10$, then $\mu'(v) \geq -1 + \frac{3}{5} + \frac{2}{5} = 0$, for otherwise there is a $\langle 5, 9, 5, 6, 10 \rangle$-star. If $\kappa_{5} = 11$, then $\mu'(v) \geq -1 + \frac{15}{22} + \frac{1}{3} \geq 0$, for otherwise there is a $\langle 5, 8, 5, 6, 11\rangle$-star. If $12 \leq \kappa_{5} \leq 17$, then $\mu'(v) \geq -1 + \frac{3}{4} + \frac{1}{4} = 0$, for otherwise there is a $\langle 5, 7, 5, 6, 17 \rangle$-star. If $\kappa_{5} \geq 18$, then $\mu'(v) \geq -1 + \frac{3}{2} \cdot \frac{18 - 6}{18} \geq 0$. 

Suppose that $v_{4}$ is a $7$-vertex. If $\kappa_{5} = 7$, then $\mu'(v) \geq -1 + 2 \cdot \frac{1}{4} + \frac{1}{2} = 0$, for otherwise there is a $\langle 5, 7, 7, 5, 11 \rangle$-star. If $\kappa_{5} = 8$, then $\mu'(v) \geq -1 + \frac{1}{4} + \frac{3}{8} + \frac{2}{5} \geq 0$, for otherwise there is a $\langle 5, 7, 8, 5, 9 \rangle$-star. If $\kappa_{5} = 9$, then $\mu'(v) \geq -1 + \frac{1}{4} + \frac{5}{12} + \frac{1}{3} = 0$, for otherwise there is a $\langle 5, 8, 5, 7, 9 \rangle$-star. If $\kappa_{5} \in \{10, 11\}$, then $\mu'(v) \geq -1 + \frac{3}{5} + 2 \cdot \frac{1}{4} \geq 0$, for otherwise there is a $\langle 5, 7, 5, 7, 11 \rangle$-star. If $\kappa_{5} \geq 12$, then $\mu'(v) \geq -1 + \frac{1}{4} + \frac{3}{4} = 0$.

Suppose that $v_{4}$ is an $8$-vertex. If $\kappa_{2} \geq 8$, then $\mu'(v) \geq -1 + 2 \cdot \frac{3}{8} + \frac{1}{4} = 0$. If $\kappa_{2} \leq 7$, then $\mu'(v) \geq -1 + \frac{3}{8} + \frac{15}{22} \geq 0$, for otherwise there is a $\langle 5, 7, 5, 8, 10 \rangle$-star. 

Suppose that $v_{4}$ is a $9$-vertex. If $\kappa_{5} = 9$, then $\mu'(v) \geq -1 + 2 \cdot \frac{5}{12} + \frac{1}{4} \geq 0$, for otherwise there is a $\langle 5, 7, 5, 9, 9 \rangle$-star. If $\kappa_{5} \geq 10$, then $\mu'(v) \geq -1 + \frac{5}{12} + \frac{3}{5} \geq 0$. 

\begin{subcase}
The $5$-vertex $v$ has precisely three $5$-neighbors $v_{2}, v_{3}$ and $v_{4}$. 
\end{subcase}

Note that the vertex $v$ is just a receiver, so we have that $\min\{\kappa_{1}, \kappa_{5}\} \neq 7, 8$. If $\min\{\kappa_{1}, \kappa_{5}\} = 6$, then $\mu'(v) \geq -1 + \frac{3}{2} \cdot \beta(18) \geq 0$, for otherwise there is a $\langle 5, 6, 6, 6, 17 \rangle$-star. If $\min\{\kappa_{1}, \kappa_{5}\} \geq 9$, then $\mu'(v) \geq 5- 6 + \frac{5}{12} + \frac{3}{5} \geq 0$, for otherwise there is a $\langle 5, 7, 5, 9, 9 \rangle$-star. 
\begin{subcase}
The $5$-vertex $v$ has precisely three $5$-neighbors $v_{2}, v_{3}$ and $v_{5}$. 
\end{subcase}

Suppose that $\kappa_{4} = 7$ and $v_{2}v_{3}$ lies in another $3$-face $uv_{2}v_{3}$. If $u$ is a $12^{+}$-vertex, then $v$ receives at least $\frac{1}{4}$ from $v_{2}$ by \ref{B1}, and then $\mu'(v) \geq -1 + \frac{1}{4} + \frac{3}{4} = 0$, for otherwise $v$ is the center of a $\langle 5, 5, 7, 5, 23 \rangle$-star. If $u$ is an $11^{-}$-vertex, then $v$ receives $\frac{1}{4}$ from $v_{3}$ by \ref{B2}, and then $\mu'(v) \geq -1 + \frac{1}{4} + \frac{3}{4} = 0$, for otherwise $v$ is the center of a $\langle 5, 5, 7, 5, 23 \rangle$-star. 

If $\min\{\kappa_{1}, \kappa_{4}\} = 8$, then $v$ receives at least $\frac{1}{8}$ from $v_{2}$ or $v_{3}$ due to \ref{B1} and \ref{B3}, which implies that $\mu'(v) \geq -1 + \frac{1}{4} + \frac{1}{8} + \frac{5}{8} \geq 0$, for otherwise there is a $\langle 5, 5, 8, 5, 15 \rangle$-star. 

If $\min\{\kappa_{1}, \kappa_{4}\} = \kappa_{1} = 9$. By the absence of $\langle 5, 5, 9, 5, 14 \rangle$-stars, we have that $\kappa_{4} \geq 15$. If $v$ is a twice-weak $5$-neighbor of $v_{1}$, then $\mu'(v) \geq -1 + \frac{1}{3} + \frac{3}{5} + \big(\frac{3}{5} - \frac{1}{2} \big) \geq 0$ by \ref{B1}; otherwise we have that $\mu'(v) \geq -1 + \big(\frac{1}{3} + \frac{1}{12}\big) + \frac{3}{5} \geq 0$. 

If $\min\{\kappa_{1}, \kappa_{4}\} = 10$, then $\mu'(v) \geq -1 + \frac{2}{5} + \frac{3}{5} = 0$, for otherwise there is a $\langle 5, 5, 10, 5, 14 \rangle$-star. If $\min\{\kappa_{1}, \kappa_{4}\} = 11$, then $\mu'(v) \geq -1 + \frac{5}{11} + \frac{4}{7} \geq 0$, for otherwise there is a $\langle 5, 5, 11, 5, 13 \rangle$-star. If $\min\{\kappa_{1}, \kappa_{4}\} \geq 12$, then $\mu'(v) \geq -1 + 2 \cdot \frac{1}{2} = 0$ by \ref{B1}. This completes the proof of \autoref{MR2}. 

\vskip 0mm \vspace{0.3cm} \noindent{\bf Acknowledgments.} This work was supported by the National Natural Science Foundation of China (xxxxxxxx) and partially supported by the Fundamental Research Funds for Universities in Henan (YQPY20140051). The authors would like to thank the anonymous referees for their valuable comments and careful reading of this paper.

\end{document}